\documentclass[11pt]{article}
\usepackage[margin=1in]{geometry}
\usepackage[utf8]{inputenc} 
\usepackage{lmodern}
\usepackage[T1]{fontenc}
\usepackage{amsmath,amssymb,subcaption,graphicx,enumitem,xcolor,tikz-cd, multicol,comment,url} 

\DeclareUnicodeCharacter{0301}{\'{}}

\usepackage{thm-restate}

\usepackage[
    eprint=true
]{biblatex}

\usepackage{xcolor}
\graphicspath{{Figures/}}

\usepackage{amsthm, tikz-cd, multicol}
\usetikzlibrary{positioning,calc}
\tikzset{main node/.style={circle,fill=black,draw,minimum size=.25cm,inner sep=0pt},}

\usepackage{hyperref}
\usepackage[capitalize,nameinlink,noabbrev,nosort]{cleveref}
\crefname{question}{Question}{Questions}

\title{Order and Gorenstein properties for lattice path matroid polytopes}
\newtheorem{theorem}{Theorem}[section]
\newtheorem{corollary}{Corollary}[theorem]
\newtheorem{lemma}[theorem]{Lemma}
\theoremstyle{definition}
\newtheorem{definition}{Definition}[section]
\theoremstyle{definition}
\newtheorem{proposition}[theorem]{Proposition}
\theoremstyle{definition}
\newtheorem{example}[theorem]{Example}
\theoremstyle{remark}
\newtheorem{remark}{Remark}
\theoremstyle{definition}
\newtheorem{question}{Question}
\theoremstyle{definition}

\DeclareMathOperator{\down}{\downarrow}

\begin{document}

\author{Carolina Benedetti, Kolja Knauer, Jer\'onimo Valencia-Porras}

\maketitle

\vspace{-0.5cm}
\begin{abstract}
    We characterize matroids whose polytopes are order polytopes as a special class of lattice path matroids, called snakes. This shows that snakes are exactly the objects lying in the intersection of the Neggers-Stanley conjecture and a conjecture of de Loera, Haws, and K{\"o}ppe. We then characterize Gorenstein lattice path matroid polytopes, yielding a new class of matroids satisfying the latter conjecture. Finally, we also show that lattice path matroids are exactly the class of positroids in which the coarsest subdivision of the matroid base polytope into translates of order polytopes is a (finest) matroidal subdivision.
       
    \end{abstract}

\section{Introduction}\label{sec:Intro}

Matroids originally arose as a combinatorial axiomatization of the concept of independence from linear algebra. They are at the core of many branches in mathematics such as graph theory, polyhedral geometry, optimization, and algebraic geometry (see~\cite{O11,ardila2010matroid,Maf81,postnikov2006total}). This manuscript is concerned with a particular class of matroids known as \emph{lattice path matroids (LPMs)} via their matroid {base} polytope, which we will call \emph{LPM polytope}. Since their introduction in~\cite{bonin2003lattice}, many different aspects of LPMs have been studied in the literature: excluded minor characterizations~\cite{Bon-10}, algebraic geometric notions~\cite{Del-12,Sch-10,Sch-11}, the Tutte polynomial~\cite{bonin2003lattice,KMR18,Mor-13}, matroid quotients~\cite{BK22,DeM-07}, and the Ehrhart polynomial~\cite{Bid-12,knauer2018lattice,FMP26}. Within the family of LPMs lies the family of \emph{snakes}. They have played a role in several of the above works in LPMs and further contain the class of minimal matroids~\cite{Ferroni_minimal}. Of particular importance is the fact that the polytope of an LPM can be split into base polytopes of snakes~\cite{chatelain2011matroid}. These are also known as \emph{shard polytopes} \cite{shardPPR23}.

We will explore LPMs in relation to other families of polytopes, such
as order polytopes introduced by Stanley~\cite{Stanley1986} and
Gorenstein polytopes as developed by Batyrev and
Nill~\cite{BatyrevNill}. Related superclasses of LPM polytopes are positroids~\cite{postnikov2006total}, polypositroids~\cite{lam2020polypositroids}, and alcoved polytopes~\cite{Lam2007}. Further, LPM polytopes are exactly the matroid polytopes that are Mirkovic--Vilonen polytopes~\cite{sanchez2023mirkovicvilonen}. Moreover, Gorenstein polytopes have been explored in relation to order polytopes ~\cite{hibi1987distributive} and matroid polytopes \cite{LM22}.
See \cref{fig:polytopeclasses} for a diagram of the related classes of polytopes mentioned in this manuscript.

\begin{figure}[ht]
    \centering
    \begin{tikzpicture}[
        dot/.style={circle, fill=black, inner sep=1.5pt},
        scale=1.5
    ]
    
        \node[dot, label=above:{\textbf{Gorenstein}}] (Gorenstein) at (-3, 3) {};
        \node[dot, label=above:{\textbf{Matroid}}] (Matroid) at (-1, 3) {};
        \node[dot, label=above:{\textbf{Alcoved}}] (Alcoved) at (1, 3) {};
        \node[dot, label=above:{\textbf{Mirković--Vilonen}}] (MV) at (3, 3) {};
    
        \node[dot, label=below left:{\textbf{\cite{LM22}}}] (33) at (-2, 2) {};
        \node[dot, label=right:{\textbf{Positroid}}] (Positroid) at (0, 2) {};
    
        \node[dot, label=below left:{\textbf{\cref{ques:PositroidGorenstein}}}] (Q1) at (-1, 1) {};
        \node[dot, label=left:{\textbf{LPM}}] (LPM) at (1, 1) {};
        \node[dot, label=right:{\textbf{Order}}] (Order) at (3, 1) {};
    
        \node[dot, label={[align=center]below left:{\textbf{\cref{thm:Gorenstein-LPM}}\\ $(\delta = 2)$}}] (Thm12_2) at (0, 0) {};
        \node[dot, label={[label distance=2pt]right:{\textbf{\cite{hibi1987distributive}}}}] (23) at (2.5, 0.5) {};
        \node[dot, label={[align=center,label distance=-12pt]below right:{\textbf{Snakes}\\(\cref{thm:snakesposets})}}] (Snakes) at (2, 0) {};
    
        \node[dot, label={[align=center]below:{\textbf{\cref{thm:Gorenstein-LPM}}\\ $(\delta\ge 3)$}}] (Thm12_3) at (1, -1) {};
    
        \newcommand{\aligncenter}[1]{\begin{tabular}{c}#1\end{tabular}}
    
        \draw (Thm12_3) -- (Thm12_2) -- (Q1) -- (33) -- (Matroid);
        \draw (Snakes) -- (LPM) -- (Positroid) -- (Alcoved);
        \draw (Thm12_2) -- (MV);
        \draw (Thm12_3) -- (Order);
    
        \draw (33) -- (Q1) -- (Thm12_2) -- (Thm12_3);
        \draw (Matroid) -- (Positroid) -- (LPM) -- (Snakes);
        \draw (Alcoved) -- (Order);
    
        \draw (Gorenstein) -- (33);
        \draw (Gorenstein) -- (23);
        \draw (Positroid) -- (Q1);
    
    \end{tikzpicture}
    \caption{Classes of polytopes related to those studied in this paper. We organize them as the Hasse diagram of the poset of containment of classes.}
    \label{fig:polytopeclasses}
\end{figure}

Our first main result shows that matroid polytopes of snakes are exactly those matroid polytopes that are also order polytopes.

\begin{restatable}{theorem}{snakesposets}
\label{thm:snakesposets}
Let $P_M$ be the polytope of a connected matroid $M$. Then the following are equivalent:
\begin{enumerate}
    \item[(i)] $M$ is a snake.
    \item[(ii)] $P_M$ is affinely equivalent to an order polytope $\mathcal{O}(X)$ of a {fence} $X$. 
    \item[(iii)] $P_M$ is affinely equivalent to an order polytope $\mathcal{O}(X)$ for some poset $X$.
    \item[(iv)] The {graph} $G(P_M)$ is isomorphic to $G(\mathcal{O}(X))$ for some order polytope $\mathcal{O}(X)$.
\end{enumerate}
\end{restatable}

In our second main result \cref{thm:Gorenstein-LPM} we characterize LPMs whose matroid polytope is Gorenstein. This theorem complements recent work from~\cite{Hibi2021,Klbl2020}. In particular, it specializes~\cite{LM22} with a very concrete description to LPMs.

\begin{restatable}{theorem}{Gorenstein}
\label{thm:Gorenstein-LPM}
Let $M=M[U,L]$ be a connected lattice path matroid on $[n]$ of rank
$k$. Then $M$ is $\delta$-Gorenstein if and only if one of the
following conditions holds.
\begin{enumerate}[leftmargin=1.5cm]
    \item[$(\delta=2)$] We have $n=2k$, the diagonal segment
    $\ell=\{(t,t):0\leq t\leq k\}$ intersects $U$ and $L$ only at
    $(0,0)$ and $(k,k)$, every concave point of $L$ is of the form
    $(i+1,i)$, and every concave point of $U$ is of the form
    $(i,i+1)$.
    \item[$(\delta\geq3)$] The matroid $M$ is
    isomorphic to $S(\delta-2,\ldots,\delta-2)$ or
    $S^*(\delta-2,\ldots,\delta-2)$.
\end{enumerate}
\end{restatable} 

It was conjectured by De Loera, Haws, and K\"oppe~\cite[Conjecture 2]{DeLoera2008} that the $h^*$-vector of any matroid polytope is unimodal. This has been proved for some uniform matroids~\cite[Theorem 3]{DeLoera2008}, sparse paving matroids of rank $2$~\cite[Theorem 1.3]{ferroni2021ehrhart}, and certain snakes~\cite{knauer2018lattice}. 
Since alcoved Gorenstein polytopes have unimodal $h^*$-vectors~\cite{bruns2007h} we have that the matroids from \cref{thm:Gorenstein-LPM} form a
new family of matroids for which the conjecture of De Loera, Haws, and K\"oppe holds true. 

On the other hand, Neggers~\cite{N78} conjectured that if $P$ is an order polytope, then its $h^*$-polynomial is real-rooted and  Stanley later even conjectured a strengthening. While the Neggers-Stanley conjecture has been disproved~\cite{B04,S07}, it remains open whether the $h^*$-vector of such $P$ has the weaker property of being unimodal, see~\cite{RW05}.
\cref{thm:snakesposets} shows that snakes are exactly the matroids lying in the intersection of the Neggers-Stanley conjecture and the conjecture of de Loera, Haws, and K\"oppe.

{As can be seen in \cref{fig:polytopeclasses}, lattice path matroids are positroids, i.e., their matroid base polytope is alcoved. The latter is equivalent to the fact that the subdivision $\Sigma_{\mathrm{box}}(M)$  induced by the hyperplanes $y_i=\sum_{j=1}^ix_j\in\mathbb Z$ yields pieces that are mapped to integer translates of order polytopes by $x\mapsto y$ as above, see~\cite{Lam2007}. In \cite{chatelain2011matroid} it is shown that if $M$ is an LPM the pieces of this subdivison are also matroidal and indeed \emph{snake polytopes}. Our third result shows that this phenomenon is restricted to LPMs.
\begin{restatable}{theorem}{OrderMatroidSubdivision}\label{prop:order-matroid-subdivision}
Let $M$ be a connected matroid. Then $M$ is a lattice path matroid if and only
if $\Sigma_{\mathrm{box}}(M)$ is matroidal. In this case,
the maximal cells of $\Sigma_{\mathrm{box}}(M)$ are snake matroid
polytopes, and $\Sigma_{\mathrm{box}}(M)$ is a finest matroid subdivision
of $P_M$.
\end{restatable}
}
{Our manuscript is structured as follows. In \cref{sec:preliminaries} we review properties of lattice path matroids and their matroid base polytopes. In \cref{sec:snakes_fences} we present relevant definitions leading to the proof of \cref{thm:snakesposets}. In \cref{sec:gorenstein} we describe further properties of LPM polytopes and characterize Gorenstein LPMs. To conclude, in \cref{sec:positroids} we discuss finest matroidal and coarsest order polytope subdivisions of matroid polytopes and show that these coincide only for the class of LPMs. 
}

A preliminary longer version of this paper has been on arXiv for a while, see~\cite{benedetti2023latticepathmatroidpolytopes}. It additionally contains a detailed analysis of the alcoved triangulation of LPM polytopes and some applications to the Ehrhart theory of LPMs of low rank. 

\section{Preliminaries}\label{sec:preliminaries}

\subsection{Lattice path matroids and  polytopes}\label{subsec:LPMsOrder}

A \textit{matroid} is a pair $M = \left(E,\mathcal{B}\right)$ where $E$ is a finite set and $\mathcal{B}$ is a non-empty collection of subsets of $E$ satisfying the following exchange axiom: 
\begin{center}
          for all $B_1,B_2 \in \mathcal{B}$, if $b_1\! \in\!  B_1\! \setminus\!  B_2$, then there is $b_2\!  \in\!  B_2\! \setminus\!  B_1$ such that $B_1\! \setminus\! \{b_1\}\!  \cup\!  \{b_2\}\!  \in\!  \mathcal{B}$. 
\end{center}
For our purposes the set $E$  is identified with $[n]:=\{1,\dots,n\}$ if $|E|=n$, and we refer to it as the \emph{ground set of $M$}. The elements in $\mathcal B$ are the \emph{bases} of $M$. It can be shown that each element in $\mathcal B$ has the same size, say $k$, and we say in this case that  the \emph{rank $r(M)$} of $M$ is $k$. The \emph{direct sum} $M_1\oplus M_2$ of two matroids $M_1=(E_1,\mathcal{B}_1), M_2=(E_2,\mathcal{B}_2)$ where $E_1\cap E_2=\emptyset$ is the matroid on $E_1\cup E_2$ whose set of bases is $\{B_1\cup B_2\mid B_1\in \mathcal{B}_1, B_2\in \mathcal{B}_2\}$. A matroid is \emph{connected}  if it cannot be written as the direct sum of two matroids.

In order to introduce the class of matroids we work with, we use the following correspondence that will be used implicitly throughout the paper. To any $k$-element subset $B\subseteq[n]$ we associate a monotone lattice path from $(0,0)$ to $(n-k,k)$ whose $i$-th step is north if and only if $i\in B$. Conversely, from any such lattice path one can recover a $k$-subset of $[n]$.

\begin{definition}\label{def:lpm} 
Let $U$ and $L$ be two lattice paths from $(0,0)$ to $(n-k,k)$ such that $U$ is above $L$. We denote by $M[U,L]$ the \emph{lattice path matroid (LPM)} that has ground set $[n]$ and its bases correspond to the set of  monotone lattice paths from $(0,0)$ to $(n-k,k)$ in the region enclosed by $U$ and $L$. We refer to this region as the \emph{diagram} of $M[U,L]$. See  Figure~\ref{fig:xmpl}.
\end{definition}

\begin{figure}[htp]
    \centering
    \begin{tikzpicture}
        \node at (0,0) {
            \begin{tikzpicture}
                 \draw[step=0.5cm, color=gray] (0,0) grid (1,1);
                 \draw[step=0.5cm, color=gray] (0.5,1) grid (1,1.5);
                 \draw[step=0.5cm, color=gray] (1,0.5) grid (1.5,1.5);
                 \draw[step=0.5cm, color=gray] (1,1.5) grid (2,2);
                \draw[line width=1pt, color = gray] (0,0) -- (0,1) -- (0.5,1) -- (0.5,1.5) -- (1,1.5) -- (1,2) -- (2,2);
                \draw[line width=1pt, color = gray] (0,0) -- (1,0) -- (1,0.5) -- (1.5,0.5) -- (1.5,1.5) -- (2,1.5)-- (2,2);

                \draw[line width = 2pt, color = black] (0,0) -- (0.5,0) -- (0.5,1) -- (1,1) -- (1,1.5) -- (1.5,1.5) -- (1.5,2) -- (2,2);
        \end{tikzpicture}
        };
        \node at (3,0) {$B = \{2,3,5,7\}$};
        \draw[->,line width = 1pt] (0.25,-0.25) to[out=-40, in=-120] (2.5,-0.35);
    \end{tikzpicture}
    \caption{Diagram of $M=M[1246,3568]$ highlighting in bold one of its bases $B$.}
    \label{fig:xmpl}
\end{figure}

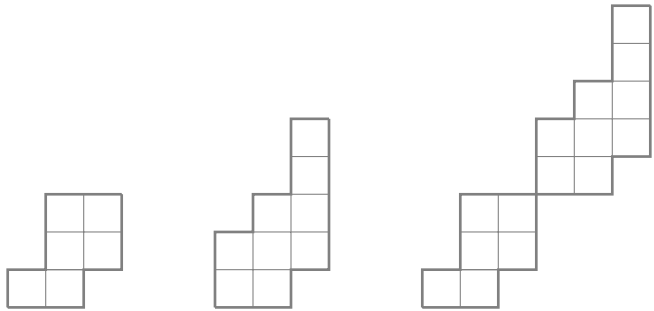
\begin{figure}[h]
    \centering
    \begin{tikzpicture}
             \draw[step=0.5cm, color=gray] (0, 0) grid (1,0.5);
             \draw[step=0.5cm, color=gray] (0.5, 0.5) grid (1.5,1.5);
            \draw[line width=1pt, color = gray] (0,0) -- (0,0.5) -- (0.5,0.5) -- (0.5,1.5) -- (1.5,1.5);
            \draw[line width=1pt, color = gray] (0,0) -- (1,0) -- (1,0.5) -- (1.5,0.5) -- (1.5,1.5);
    \end{tikzpicture}
    \hspace{1cm}
    \begin{tikzpicture}
             \draw[step=0.5cm, color=gray] (0,0) grid (1,1);
             \draw[step=0.5cm, color=gray] (0.5,1) grid (1,1.5);
             \draw[step=0.5cm, color=gray] (1,0.5) grid (1.5,2.5);
            \draw[line width=1pt, color = gray] (0,0) -- (0,1) -- (0.5,1) -- (0.5,1.5) -- (1,1.5) -- (1,2.5) -- (1.5,2.5);
            \draw[line width=1pt, color = gray] (0,0) -- (1,0) -- (1,0.5) -- (1.5,0.5) -- (1.5,2.5);
    \end{tikzpicture}
    \hspace{1cm}
    \begin{tikzpicture}
            \draw[step=0.5cm, color=gray] (0, 0) grid (1,0.5);
            \draw[step=0.5cm, color=gray] (0.5, 0.5) grid (1.5,1.5);
             \draw[step=0.5cm, color=gray] (1.5,1.5) grid (2.5,2.5);
             \draw[step=0.5cm, color=gray] (2.5, 2) grid (3,4);
             \draw[step=0.5cm, color=gray] (2,2.5) grid (2.5,3);
            \draw[line width=1pt, color = gray] (0,0) -- (0,0.5) -- (0.5,0.5) -- (0.5,1.5) -- (1.5,1.5);
            \draw[line width=1pt, color = gray] (0,0) -- (1,0) -- (1,0.5) -- (1.5,0.5) -- (1.5,1.5);
            \draw[line width=1pt, color=gray] (1.5,1.5) -- (1.5,2.5) -- (2,2.5) -- (2,3) -- (2.5,3) -- (2.5,4) -- (3,4);
            \draw[line width=1pt, color=gray] (1.5,1.5) -- (2.5,1.5) -- (2.5,2) -- (3,2) -- (3,4);
    \end{tikzpicture}
    \caption{Diagrams of two LPMs and their direct sum.}
    \label{fig:LPM-direct-sum}
\end{figure}

Given two lattice paths $U,L$ from $(0,0)$ to $(n-k,k)$, 
the rank of $M[U,L]$ is $k$, i.e., the height of the diagram. Further, $M[U,L]$ is {connected} if and only if $U$ and $L$ only intersect in $(0,0)$ and $(n-k,k)$. Indeed, note that the direct sum of two LPMs corresponds to the concatenation of their corresponding diagrams, as shown in Figure~\ref{fig:LPM-direct-sum}.

\begin{definition}\label{def:snakes}
A lattice path matroid $S$ is a \emph{snake} if $S$ is connected and its diagram does not contain a $2\times 2$ square. 
\end{definition}
 
We now give a more direct description for snakes. Let $\alpha=(\alpha_1,\dots,\alpha_s)$ be a tuple of positive integers such that $\alpha_1+\cdots+\alpha_s=n$. That is, $\alpha$ is an \emph{integer (strong) composition of $n$}, denoted $\alpha\models n$. Identify each $\alpha_i$ with a rectangle $R_i$ of size $1\times (\alpha_i+1)$, called a \emph{vertical strip}. Such $R_i$ gives rise to a \emph{horizontal strip} $R_i^*$ which is the rectangle $(\alpha_i+1)\times 1$. Notice that $R_i$ and $R_i^*$ contain $\alpha_i+1$ unit squares.
Given a composition $\alpha=(\alpha_1,\dots,\alpha_s)$ let $S(\alpha)$ be the snake whose diagram is obtained as follows. 
From $(0,0)$ draw $R_1$, then identify its last square with the first square of $R_2^*$, then the last square of $R_2^*$ with the first square of $R_3$, and so on. On the other hand, the composition $\alpha=(\alpha_1,\dots,\alpha_s)$ gives rise to the snake $S^*(\alpha)$ obtained by reflecting $S(\alpha)$ across the line $y=x$. More precisely, $S^*(\alpha)$ is obtained as follows: from $(0,0)$ draw $R_1^*$, then identify its last square with the first square of $R_2$, then the last square of $R_2$ is identified with the first square of $R_3^*$, and so on. Notice that if $\alpha\models n$ then $S(\alpha)$, and hence $S^*(\alpha)$, are LPMs on the set $[n+2]$. See Figure \ref{fig:snake-and-poset} for an example of a snake $S(\alpha).$

A matroid $M=([n],\mathcal B)$ can be characterized geometrically via its \emph{matroid base polytope}, which we denote by $P_M$. This is a lattice polytope and its vertices are given by $\{e_B: B\in\mathcal B\}$ where $e_B:= \sum_{i\in B}{e}_i$ and $\{e_1,\dots,e_n\}$ is the canonical basis of $\mathbb R^n$. Moreover, the edges of the polytope encode the exchange axiom described before. In view of this, we sometimes think of $M$ as being $P_M$, and vice versa (see \cite{GGMS}).
It is known that if $M$ is an LPM then its matroid polytope $P_M$ is an alcoved polytope \cite{Lam2007, knauer2018lattice}.
Further, it is known that $P_M$ can be decomposed into subpolytopes, where each subpolytope in the decomposition is the matroid polytope of a {snake}~\cite{chatelain2011matroid}. 

Our results shed light on the base polytopes of LPMs and are motivated by understanding their $h^*$-vectors. Thus, even if we do not use its definition explicitly, let us define it.

Given a {lattice polytope} $P$, let $tP$ be its dilation by an integer factor of $t$. The number of integer points in $tP$ is known to be the evaluation $L_P(t)$ of the \emph{Ehrhart polynomial} of $P$~\cite{Ehrhart}. The generating series of the sequence $\{ L_P(t)\}_{t\geq 0}$, known as the Ehrhart series of $P$, can be written as a rational function 
$$
\mathrm{Ehr}_P(z):=\displaystyle\sum_{t\geq 0}L_P(t)z^t=\dfrac{h^*_0+h^*_1z+\cdots+h^*_mz^m}{(1-z)^{d+1}}.
$$
The polynomial $h^*(P):=h^*_0+h^*_1z+\cdots+h^*_mz^m$ is known as the \emph{$h^*$-polynomial of $P$}, the tuple $(h_0^*,h_1^*,\dots,h_m^*)$ is the \emph{$h^*$-vector of $P$}, and the degree of $h^*$ satisfies $m \leq \dim(P)$. 

Stanley's non-negativity theorem \cite{Stanley-NonNegThm} guarantees that if $P$ is a lattice polytope then its $h^*$-vector has integer and nonnegative coordinates. Thus, it is natural to wonder if the coefficients of $h^*$ have a combinatorial interpretation or satisfy any natural properties of integer sequences such as being unimodal.

\section{Snakes and fences}\label{sec:snakes_fences}
An ideal of a poset $X$ is a subset $I$ such that if $x\in I$ and $y\leq x$ then $y\in I$. Following~\cite{Stanley1986}, the \emph{order-polytope} $\mathcal{O}(X)$ is the convex hull of all the $e_I$ for $I$ an ideal of $X$.\footnote{Note that Stanley defines $\mathcal{O}(X)$ as convex hull of filters, but this is affinely equivalent.}
The main result in this section is that matroid polytopes that are affinely equivalent to order polytopes are matroid polytopes of snakes and these are the order polytopes of fences.

\begin{definition}\label{def:fence_posets}
    Let $\alpha=(\alpha_1,\dots,\alpha_s)\models n$. 
    The \emph{fence} of $\alpha$, denoted $F(\alpha)$ is the poset on $n+1$ elements $p_1,\ldots,p_{n+1}$ with the covering relations: 
    \begin{equation*}
    p_1\prec p_2 \prec \cdots\prec p_{\alpha_1+1}\succ p_{\alpha_1+2}\succ \cdots\succ p_{\alpha_1+\alpha_2+1}\prec p_{\alpha_1+\alpha_2+2}\prec\cdots\prec p_{\alpha_1+\alpha_2+\alpha_3+1}\succ \cdots.
    \end{equation*}
\end{definition}

\begin{figure}[ht]
    \centering
    \begin{tikzpicture}
            \draw[step=0.5cm, color=gray] (0, 0) grid (0.5,1);
            \draw[step=0.5cm, color=gray] (0.5, 0.5) grid (1,1.5);
            \draw[step=0.5cm, color=gray] (1,1) grid (1.5,2.5);
            
            \draw[line width=2pt, color = black] (0,0.5) -- (0,1) -- (0.5,1) -- (0.5,1.5) -- (1,1.5) -- (1,2.5);
            \draw[line width=1pt, color=gray] (0,0) -- (0.5,0) -- (0.5,0.5) -- (1,0.5) -- (1,1) -- (1.5,1) -- (1.5,2.5);
            \draw[line width=1pt, color = gray] (0,0) -- (0,0.5);
            \draw[line width=1pt, color = gray] (1,2.5) -- (1.5,2.5);
    \end{tikzpicture}
    \hspace{2cm}
    \begin{tikzpicture}[rotate = -45]
            \draw[line width=1pt] (0.5,0) -- (0.5,0.5) -- (1,0.5) -- (1,1) -- (1.5,1) -- (1.5,2); 
            
            \filldraw [color = black] (0.5,0) circle (2pt);
            \filldraw [color = black] (0.5,0.5) circle (2pt);
            \filldraw [color = black] (1,0.5) circle (2pt);
            \filldraw [color = black] (1,1) circle (2pt);
            \filldraw [color = black] (1.5,1) circle (2pt);
            \filldraw [color = black] (1.5,1.5) circle (2pt);
            \filldraw [color = black] (1.5,2) circle (2pt);
    \end{tikzpicture}
    \caption{The snake $S(1,1,1,1,2) = M[12467,24678]$ and its associated fence $F(1,1,1,1,2)$.}
    \label{fig:snake-and-poset}
\end{figure}
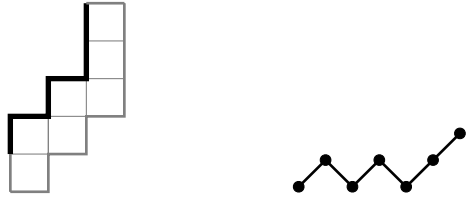

\begin{remark}
For our purposes, the poset $F^*(\alpha)$ dual to $F(\alpha)$ is also considered a fence.
\end{remark}
Fences are a natural class of posets that appear in the study of cluster algebras, quiver representations, and other areas of enumerative combinatorics, see~\cite{MSS21} for an overview.  In~\cite{knauer2018lattice}, fences are called \emph{zig-zag-chain posets}. Our choice of notations for snakes and fences reflects the correspondence between them. More concretely, one way of constructing the Hasse diagram of the fence poset from the diagram of the snake is to delete the first and last steps in the upper (or lower) path of the diagram of $S(\alpha)$ and rotate this path $45$ degrees clockwise. We describe this construction through an example in Figure~\ref{fig:snake-and-poset}.

In~\cite[Theorem 4.7]{knauer2018lattice} it is shown that the matroid base polytope $P_{S(\alpha)}$ of the snake $S(\alpha)$ and the order polytope $\mathcal{O}(X)$ of the fence $F(\alpha)$ are affinely equivalent. Moreover, this result extends to direct sums of snakes and disjoint unions of fences.

 Before stating our main result, further notation is in order. Given a polytope $P$ we denote by $G(P)$ its 1-skeleton. In other words, $G(P)$ is the graph of the polytope $P$. 

\begin{lemma}\cite[Lemma 1.1a]{HIBI2017991}\label{lemma:edges_order_polytope}
Let $I$ and $J$ be ideals of a poset $X$, with $I\neq J$. Then $\{e_I,e_J\}$ is an edge in $G(\mathcal{O}(X))$ if and only if $I\subset J$ and $J\setminus I$ induces a connected subposet of $X$.
\end{lemma}

Given a poset $X$ and $x\in X$ we denote by $\down x$ the principal ideal generated by $x$. That is, $\down x:=\{y\in X : y\leq x \}$. 
As a consequence of Lemma~\ref{lemma:edges_order_polytope} the reader is encouraged to verify that the following claims hold:

\begin{itemize}
    \item[$\circ$] If there is no containment relation between $I$ and $J$, then vertices $e_I,e_J$ are not adjacent.
    \item[$\circ$] If $J=\down x$ for some $x\in X$ and $I\subset J$, then $\{e_I,e_J\}$ is an edge of $\mathcal{O}(X)$ since for every $y\in J\setminus I$ there is a path from $y$ to $x$, and thus $J\setminus I$ is connected. More generally, if $J=\down x\cup I$ and $x\notin I$, then $\{e_I,e_J\}$ is an edge of $\mathcal{O}(X)$.
\end{itemize}

For a graph $G$, let $V_G$ denote the set of its vertices, and for $x,y\in V_G$, let $d(x,y)$ be the \emph{distance between $x$ and $y$ on $G$}, defined as the length of the shortest path between the vertices on the graph. Moreover, denote by $I(x,y)$ the \emph{interval} in $G$ corresponding to the subgraph induced by all the vertices that lie on a shortest path. In other words, $I(x,y)=\{z\in V_G\mid d(x,z)+d(z,y)=d(x,y)\}$.

\begin{lemma}\cite[Lemma 1.4, Lemma 1.6]{MAURER1973216}\label{lem:IC-PC}
    Let $M$ be a matroid and $P_M$ its matroid polytope. The following conditions hold:
    \begin{itemize}
        \item \textbf{Interval Condition (IC)}: If $x,y$ are vertices of $G(P_M)$ with $d(x,y) = 2$, then $I(x,y)$ is (isomorphic to) the graph of a square, a square pyramid or an octahedron.
        \item \textbf{Positioning Condition (PC)}: If $a,b,c,d$ are vertices that induce a square in $G(P_M)$, then for any vertex  $x$ it holds that $d(x,a) + d(x,c) = d(x,b) + d(x,d)$.
    \end{itemize}
\end{lemma}

Now we are in position to state our main result in this section.

\snakesposets*

\begin{proof} 
As mentioned before, conditions ($i$) and ($ii$) are proven to be equivalent in~\cite[Theorem 4.7]{knauer2018lattice}. Also, it is clear that ($ii$)$\implies$ ($iii$) and ($iii$)$\implies$ ($iv$). 
Now we will show that ($iv$)$\implies$($ii$) and this will complete the proof.

Let $M$ be a matroid with polytope $P_M$ such that the graph $G(P_M)$ is isomorphic to $G(\mathcal{O}(X))$ for some poset $X$. Given a poset $X$ we denote by $G_X$ the undirected graph obtained from the Hasse diagram of $X$, ignoring orientation of the edges. The graph $G_X$ is known as the \emph{cover graph of $X$}. Thus, a poset $X$ is a fence if its cover graph $G_X$ is a path. Thus, it remains to show that the cover graph $G_X$ is a path.
Indeed, suppose that this has been established. Then $X$ is a fence,
and by \cite[Theorem~4.7]{knauer2018lattice} the order polytope
$\mathcal O(X)$ is affinely equivalent to the base polytope $P_S$ of
a snake $S$. Hence
$G(P_M)\cong G(\mathcal O(X))\cong G(P_S)$.

By~\cite{HolzmannNortonTobey} or \cite[Theorem~13]{PinedaVillavicencioSchroeter}, an abstract basis
exchange graph determines its matroid up to relabelling, duality of
connected components, and the addition or deletion of loops and
coloops. Each of these operations preserves the affine equivalence
class of the base polytope: relabelling permutes the coordinates,
duality on a component replaces $x_i$ by $1-x_i$ on that component,
and loops and coloops correspond to fixed coordinates. Therefore
$P_M$ is affinely equivalent to $P_S$, and hence to
$\mathcal O(X)$. This proves~(ii).

Thus, we will show that the cover graph of $X$, $G_X$, is such that the maximum degree among all of its vertices is $2$, and that $G_X$ is acyclic and connected. Since isomorphic basis graphs determine matroid polytopes up to affine equivalence this suffices. The proof goes by contradiction, considering different cases which we illustrate in Figure~\ref{fig:snakeorder}.

\begin{figure}[htp]
    \centering
    \includegraphics[width=\textwidth]{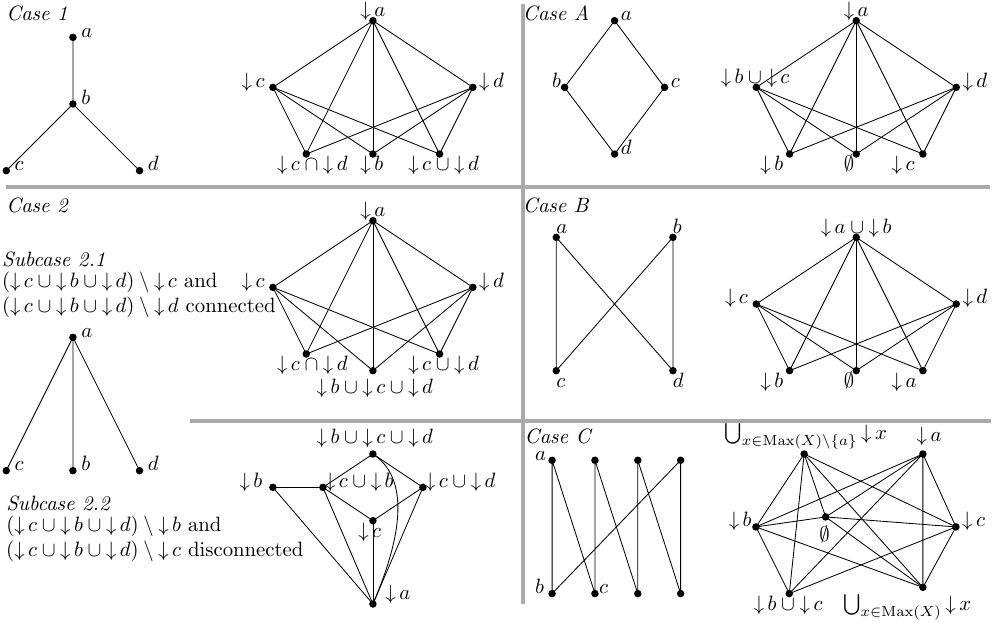}
    \caption{The cases in the proof of Theorem~\ref{thm:snakesposets}.}\label{fig:snakeorder}
\end{figure}

As a first observation one can see that if $G_X$ has several connected components corresponding to different posets, then $\mathcal{O}(X)=P_M$ would be a product of their order polytopes. However, this would contradict that $M$ is connected. 

Now, suppose that $G_X$ has a vertex of degree larger than $2$. There are essentially two ways this can happen:

\medskip

\noindent\emph{Case 1}: There are $c,d\prec b\prec a$ in $X$ and $c,d$ are incomparable. In this case consider the four principal ideals $\down a$, $\down  b$, $\down  c$, $\down  d$, along with $\down c \cup \down d$ and $\down c \cap \down d$. Identifying the ideal $\down  x$ with its corresponding vertex in $G(\mathcal{O}(X))$ we can see via Lemma \ref{lemma:edges_order_polytope} that $d(\down  c,\down  d)=2$ and all these ideals belong to the interval $[\down  c,\down  d]$. However, the vertex $\down a$ has degree at least $5$. Hence, the interval $[\down  c,\down  d]$ does not induce a square, square pyramid, or octahedron. This contradicts the IC condition from~\cref{lem:IC-PC} and thus $G(\mathcal{O}(X))$ cannot be isomorphic to $G(P_M)$.

\medskip

\noindent \emph{Case 2}: 
There are $b,c,d \prec a$ in $X$ and $b,c,d$ are incomparable.

\medskip

\noindent \emph{Subcase 2.1}: $(\down c \cup\down b \cup \down d) \setminus \down c$ and $(\down c \cup\down b \cup \down d) \setminus \down d$  both induce connected subposets of $X$.

Using Lemma \ref{lemma:edges_order_polytope} we can see that $d(\down  c,\down  d)=2$ and the interval $[\down  c,\down  d]$ contains $\down a$, $\down c \cup\down b \cup \down d$, $\down  c$, $\down  d$, and further $\down c \cup \down d$, and $\down c \cap \down d$. As in the previous case the vertex $\down a$ has degree at least $5$. Hence, the interval $[\down  c,\down  d]$ does not induce a square, square pyramid, or octahedron. Again, this contradicts the IC condition and thus $G(\mathcal{O}(X))\ncong G(P_M)$.

\medskip

After possible relabeling it only remains to prove one further subcase:
\medskip

\noindent \emph{Subcase 2.2}: 
     If $(\down c \cup\down b \cup \down d) \setminus \down b$ and $(\down c \cup\down b \cup \down d) \setminus \down c$ are disconnected then $\down c, \down c \cup\down b, \down c \cup\down b \cup \down d, \down c \cup\down d$ induce a square, using in particular the assumption that $\down c \cup\down b \cup \down d\setminus \down c$ is disconnected. Now, distances from  $\down b$ yield
    \begin{itemize}
        \item $d(\down b, \down c) = 2$ since they are non-adjacent but can be connected through $\down a$
        \item $d(\down b , \down c \cup\down b) = 1$ as they form an edge,
        \item $d(\down b, \down c \cup\down d) = 2$  since they are non-adjacent but can be connected through $\down a$, 
        \item $d(\down b, \down c \cup\down b \cup \down d) = 2$ since they are non-adjacent by assumption but can be connected through $\down a$.  
    \end{itemize} 
Altogether, this contradicts the PC condition from Lemma~\ref{lem:IC-PC}.

Note that if $X^*$ denotes the dual poset of $X$ then $G(\mathcal{O}(X))\cong G(\mathcal{O}(X^*))$. 
Hence, the maximum degree of $G_X$ is at most $2$.

It now only remains to show that $G_X$ contains no cycles. By contradiction suppose $G_X$ has a cycle, then since the maximum degree of $G_X$ is $2$, $G_X$ must be a cycle $C$. It holds that $X$ has the same number of minima and maxima and thus we distinguish three cases:

\medskip

\noindent \emph{Case A}: $X$ has exactly one maximum:
    Let $a$ be the maximum, $d$ the minimum and $b,c$ two incomparable elements such that $d \leq b,c \leq a$. Then $d(\down d, \down b\cup \down c)=2$ and the interval $I(\down d, \down b\cup \down c)$ contains the vertex $\down a$ whose degree is at least $5$ since it is adjacent to $\down d$, $\down b$, $\down c$,  $\down a\cup \down b$ and $\emptyset$. This, and the fact that all of these vertices belong to the interval $I(\down d, \down b\cup \down c)$ contradict IC from \cref{lem:IC-PC}. 

\medskip

\noindent \emph{Case B}: $X$ has exactly two maxima:
Let $a,b$ be the maxima, and $c,d$ the minima of $X$. Then $d(\down c, \down d)=2$ and  $I(\down c, \down d)$ contains vertex $\down a\cup \down b$ whose degree is at least $5$. Indeed, $\down a\cup \down b$ is adjacent to $\down a$, $\down b$, $\down c$, $\down d$, $\emptyset$.  Again, this violates IC from Lemma \ref{lem:IC-PC}.

\medskip
\noindent \emph{Case C}: $X$ has at least three maxima:
Denote by $\mathrm{Max}(X)$ and $\mathrm{Min}(X)$ the set of maxima and minima of $X$, respectively. Let $a\in \mathrm{Max}(X)$ be a maximum and $b,c\in \mathrm{Min}(X)$ such that $b,c\leq a$. Now,  $d(\down b, \down c)=2$ and the interval $I(\down b, \down c)$ contains vertices $\bigcup_{x\in \mathrm{Max}(X)}\down x$, $\bigcup_{x\in \mathrm{Max\setminus\{a\}}(X)}\down x$, $\down a$, $\down b$, $\down b\cup \down c$, $\down c$, $\emptyset$. This  contradicts IC from Lemma \ref{lem:IC-PC}. 
Hence, $G_X$ is acyclic and the result follows.

\end{proof}

\section{Gorenstein LPMs}\label{sec:gorenstein}

We first fix the lattice conventions used in this section. Let
$P\subseteq\mathbb R^n$ be a lattice polytope and $V_P$  the linear space parallel to $\operatorname{aff}(P)$, and 
$N_P:=V_P\cap\mathbb Z^n$. Thus, $P$ is regarded as a
full-dimensional polytope in its affine span, endowed with the affine
lattice $\operatorname{aff}(P)\cap\mathbb Z^n$. If
$z\in\operatorname{aff}(P)\cap\mathbb Z^n$, then $P-z\subseteq V_P$,
and $N_P$ is the lattice induced on $V_P$. We write
$N_P^\vee:=\operatorname{Hom}(N_P,\mathbb Z)$ for the dual lattice.

Suppose that $0\in\operatorname{relint}(P)$. The \emph{polar} of $P$ in
$V_P$ is
$P^\vee:=\{u\in V_P^*:u(x)\geq-1\text{ for every }x\in P\}$. We call
$P$ \emph{reflexive} if $P^\vee$ is a lattice polytope with respect to
$N_P^\vee$. More generally, for a positive integer $\delta$, we say
that $P$ is \emph{$\delta$-Gorenstein} if there exists
$z\in\operatorname{relint}(\delta P)\cap\mathbb Z^n$ such that
$\delta P-z$ is reflexive with respect to the lattice $N_P$. We call
$P$ \emph{Gorenstein} if it is $\delta$-Gorenstein for some $\delta$. This
relative-lattice formulation, including the non-full-dimensional
case, agrees with the usual definition, see
\cite[Section~2.3]{LM22}.

For a matroid $M$, we say that $M$ is \emph{($\delta$-)Gorenstein} if its base
polytope $P_M$ is ($\delta$-)Gorenstein. Notice that if $M$ is connected
of rank $k$ on $[n]$, then
$\operatorname{aff}(P_M)=\{x\in\mathbb R^n:
\sum_{i=1}^n x_i=k\}$ and the lattice parallel to this affine space is
$\{x\in\mathbb Z^n:\sum_{i=1}^n x_i=0\}$. Thus, all notions of
duality, primitivity, and lattice distance below are taken in this
induced lattice, rather than in the ambient lattice $\mathbb Z^n$.

In this section we will analyze LPM polytopes in order to provide a combinatorial characterization of those that are {Gorenstein}. The road map to achieve this characterization will require us to provide an understanding of minors of LPMs, as well as a hyperplane description of $P_M$. 

\subsection{Minors of LPMs}
Deletion and contraction are well known operations that can be performed on any matroid to produce a smaller matroid. These operations may be defined as follows.

Let $M=([n],\mathcal B)$ be a matroid of rank $k$ and recall that $i\in[n]$ is a \emph{coloop} if it is in every basis, and a \emph{loop} if it is in no basis.
Let $i\in[n]$ such that $i$ is not a coloop. The \emph{deletion} of $i$ from $M$ is the matroid $M\setminus i=([n]-\{i\},\mathcal B')$ where $\mathcal B'=\{B\in\mathcal B\; :\; i\notin B\}$. Thus in this case $M\setminus i$ has rank $k$.
 Let $i\in[n]$ such that $i$ is not a loop. The \emph{contraction} of $i$ from $M$ is the matroid $M/ i=([n]-\{i\},\mathcal B'')$ where $\mathcal B''=\{B-\{i\}\; :\; B\in\mathcal B, i\in B\}$. Thus in this case $M/i$ has rank $k-1$.
If $i$ is a coloop (loop) then $M\setminus i$ (respectively $M/i$) has bases as in $\mathcal B''$ (respectively $\mathcal B'$). A \emph{minor} of a matroid $M$ is a matroid obtained from $M$ by successively deleting and contracting elements.

The family of LPMs is known to be closed under minors. That is, given an LPM $M$ it holds that any minor of $M$ is again an LPM (see ~\cite{bonin2003lattice,BONIN2006701} for details). Let us give some intuition on what deleting or contracting an element in $M$, a given LPM, looks like. This intuition will be useful to understand facets of the polytope $P_M.$

Recall that the diagram of $M[U,L]$ is completely determined by the lattice paths  $U$ and $L$. Both of these paths have labels in $[n]$, and if the rank of $M$ is $k$ then $k$ of these labels are north steps in both paths. First, observe that $i\in[n]$ is a coloop if it corresponds to a north segment that is in both $U$ and $L$ and it is a loop if it is an east segment that is in both $U$ and $L$. In both cases contraction and deletion of $i$ just corresponds to contracting the segment of the diagram.

\begin{figure}[h]
    \centering
    \begin{tikzpicture}
             \draw[step=0.5cm, color=gray] (0, 0) grid (1.5,1);
             \draw[step=0.5cm, color=gray] (0.5,1) grid (2,2.5);
            \draw[line width=1pt, color=gray] (0,0) -- (0,1) -- (0.5,1) -- (0.5,2.5) -- (2,2.5);
            \draw[line width=1pt, color=gray] (0,0) -- (1.5,0) -- (1.5,1) -- (2,1) -- (2,2.5);
            
            \draw[line width=2pt, color = red] (0.5,2) -- (0.5,1.5);
            \draw[line width=2pt, color = red] (1.5,1) -- (1.5,0.5);
        
    \end{tikzpicture}
    \hspace{1cm}
    \begin{tikzpicture}
             \draw[step=0.5cm, color=gray] (0, 0) grid (1.5,1);
             \draw[step=0.5cm, color=gray] (0.5,1) grid (2,2.5);
            \draw[line width=1pt, color=gray] (0,0) -- (0,1) -- (0.5,1) -- (0.5,1.5);
            \draw[line width=1pt, color=gray] (1,2.5) -- (2,2.5);
            \draw[line width=1pt, color=gray] (0,0) -- (1,0);
            \draw[line width=1pt, color=gray] (1.5,1) -- (2,1) -- (2,2.5);
            \draw[line width=2pt, dashed,color = red] (0.5,1.5) -- (0.5,2.5) -- (1,2.5);
            \draw[line width=2pt, dashed,color = red] (1,0) -- (1.5,0) -- (1.5,1);
    \end{tikzpicture}
    \hspace{1cm}
    \begin{tikzpicture}
             \draw[step=0.5cm, color=gray] (0, 0) grid (1,1);
             \draw[step=0.5cm, color=gray] (0.5,1) grid (1.5,2.5);
            \draw[line width=1pt, color=gray] (0,0) -- (0,1) -- (0.5,1) -- (0.5,1.5);
            \draw[line width=1pt, color=gray] (0.5,2.5) -- (1.5,2.5);
            \draw[line width=1pt, color=gray] (0,0) -- (1,0);
            \draw[line width=1pt, color=gray] (1,1) -- (1.5,1) -- (1.5,2.5);
            \draw[line width=2pt,color = red] (0.5,1.5) -- (0.5,2.5);
            \draw[line width=2pt,color = red] (1,0) -- (1,1);
    \end{tikzpicture}
    \caption{Deletion of the element $5$ in M[12456,45789].}
    \label{fig:LPM-deletion}
\end{figure}
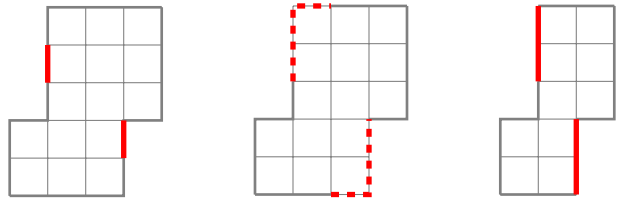

Suppose we want to obtain the diagram of the LPM obtained from $M$ by deleting the element $i$ that is not a coloop. That is, $M\setminus i$, for some $i\in[n]$. 
Start with paths $U$ and $L$ drawn bold. Then in order to obtain $M\setminus i$, a portion of $U$ and $L$ will turn dashed. The rightmost bold portions of the lattice paths $U$ and $L$ will be shifted to the right, giving us the diagram of $M\setminus i$. More precisely,
if the $i$-th step in $U$ is vertical, then dash the steps with labels $\{i,i+1,\dots,i+j \}$ where $i+j$ is the first horizontal step from $i$. If the $i$-th step in $L$ is vertical, then dash the steps $\{i,i-1,\dots,i-j \}$ where $i-j$ is horizontal and $\{i,i-1,\dots,i-j+1 \}$ are vertical in $L$.
If $i$ happens to be horizontal in $U$ (or $L$), we dash it only. Then shift to the left the right-most bold portions of $U$ and $L$, as shown in Figure \ref{fig:LPM-deletion}.

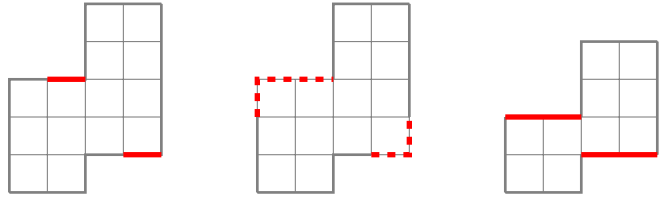
\begin{figure}[h]
    \centering
    \begin{tikzpicture}
             \draw[step=0.5cm, color=gray] (0, 0) grid (1,1.5);
             \draw[step=0.5cm, color=gray] (1,0.5) grid (2,2.5);
            \draw[line width=1pt, color=gray] (0,0) -- (0,1.5) -- (1,1.5) -- (1,2.5) -- (2,2.5);
            \draw[line width=1pt, color=gray] (0,0) -- (1,0) -- (1,0.5) -- (2,0.5) -- (2,2.5);
            
            \draw[line width=2pt, color = red] (0.5,1.5) -- (1,1.5);
            \draw[line width=2pt, color = red] (2,0.5) -- (1.5,0.5);
        
    \end{tikzpicture}
    \hspace{1cm}
    \begin{tikzpicture}
             \draw[step=0.5cm, color=gray] (0, 0) grid (1,1.5);
             \draw[step=0.5cm, color=gray] (1,0.5) grid (2,2.5);
            \draw[line width=1pt, color=gray] (0,0) -- (0,1);
            \draw[line width=1pt, color=gray] (1,1.5) -- (1,2.5) -- (2,2.5);
            \draw[line width=1pt, color=gray] (0,0) -- (1,0) -- (1,0.5) -- (1.5,0.5);
            \draw[line width=1pt, color=gray] (2,1) -- (2,2.5);
            \draw[line width=2pt, dashed, color = red] (0,1) -- (0,1.5) -- (1,1.5);
            \draw[line width=2pt, dashed, color=red] (1.5,0.5) -- (2,0.5) -- (2,1);
    \end{tikzpicture}
    \hspace{1cm}
    \begin{tikzpicture}
              \draw[step=0.5cm, color=gray] (0, 0) grid (1,1);
             \draw[step=0.5cm, color=gray] (1,0.5) grid (2,2);
            \draw[line width=1pt, color=gray] (0,0) -- (0,1);
            \draw[line width=1pt, color=gray] (1,1) -- (1,2) -- (2,2);
            \draw[line width=1pt, color=gray] (0,0) -- (1,0) -- (1,0.5);
            \draw[line width=1pt, color=gray] (2,0.5) -- (2,2);
            
            \draw[line width=2pt,color = red] (0,1) -- (1,1);
            \draw[line width=2pt,color = red] (1,0.5) -- (2,0.5);
    \end{tikzpicture}
    \caption{Contraction of the element $5$ in $M[12367,36789]$.}
    \label{fig:LPM-contraction}
\end{figure}

Now, suppose we want to contract a non-loop element $i$, and see how the contraction $M/i$ can be thought of. 
If $i$ is horizontal in $U$, then dash the steps with labels $\{i,i-1,\dots,i-j \}$ where $i-j$ is vertical and $\{i,i-1,\dots,i-j+1 \}$ are horizontal in $U$.
if $i$ is horizontal in $L$, then the dash steps $\{i,i+1,\dots,i+j \}$ where $i+j$ is the first vertical step from $i$ in $L$. 
If $i$ happens to be vertical in $U$ (or $L$), we dash it only. Then shift down the upper portions of $U$ and $L$. See Figure \ref{fig:LPM-contraction} for an illustrative example.  

\subsection{Face structure of LPM polytopes}
It is known that given a matroid $M$ and its polytope $P_M$, every face of $P_M$ is (the polytope of) a matroid (see \cite{GGMS}).
In this section we aim to describe these facets precisely when $M$ is an LPM. In particular,  we will see that the faces of $P_M$ are again LPMs. 
Given an LPM over $[n]$ of rank $k$, say $M=M[U,L]$, we have previously described $U=\{u_1,\dots,u_k\}$ by the set of labels of the $k$ north steps that define the lattice path $U$, and similar for $L$. We can also define the same lattice path $U$ as a $0/1$-vector $(U_1,U_2,\dots,U_n)$ where $U_i=1$ if and only if the $i$-th step of $U$ is north, then exactly $k$ of these $U_i$'s are equal to 1, and similar for $L$. 

With this set up in mind we establish the following result which appears in \cite[Theorem 3.3]{knauer2018lattice}, and which gives us a description of defining hyperplanes for $P_M$.

\begin{theorem}[{\cite[Theorem 3.3]{knauer2018lattice}}]\label{H-description-LPM-polytopes}
Let $M=M[U,L]$ be an $LPM$ of rank $k$ over $[n]$ such that $U = (U_1\ldots,U_{n})$ and $L = (L_1,\ldots,L_n)$. Then 
$$P_M = \left\{ x\in\mathbb{R}^n \;\;\; \Big| \;\;\; 0\leq x_i \leq 1 \;\; \text{and} \;\; \sum_{j=1}^{i} L_j \leq\sum_{j=1}^{i} x_j \leq \sum_{j=1}^{i} U_j \;\; \text{for all }\;i\in[n] \right\}.$$
\end{theorem}

We will make use of Theorem \ref{H-description-LPM-polytopes} along with the description of minors of LPMs discussed before to describe the facet-defining hyperplanes of LPMs. Since this result is already present in the literature~\cite[Proposition 14]{An17}, but in a slightly different set-up we only present the ideas here. Our characterization will be based on certain lattice points over the paths $U$ and $L$. Notice that the path $U$ (as well as $L$) contains $n$ lattice points, excluding $(0,0)$. The $i$-th lattice point in $U$ is the point $p_{i,U}$ with coordinates $p_{i,U}=(U_{0,i},U_{1,i})$ where 

$$U_{0,i}:=\displaystyle\sum_{j\leq i, U_j=0}1 \;\;\;\text{and}\;\;\; U_{1,i}:=\displaystyle\sum_{j\leq i, U_j=1}1.$$ 

Similarly, we obtain the $i$-th lattice point $p_{i,L}$ of the lattice path $L$ as $p_{i,L}=(L_{0,i},L_{1,i})$. A \emph{point on the boundary of $M$} is a lattice point $p$ lying either on $U$ or $L$. Any point $p$ on the boundary of $M$ defines four quadrants just by translating the origin to $p$.
 We call a point $p$ on the boundary of the diagram of $M=M[U,L]$ \emph{concave} if exactly one of the four quadrants at $p$ has empty intersection with the diagram of $M$. In Figure \ref{fig:concave} the point $p=(1,2)$ of $M[12467,45678]$ is concave. 
 
 Also, if $M$ is a matroid on $[n]$ and $M=M_1\oplus\cdots\oplus M_r$ is its decomposition into connected components, then $\dim P_M=n-r$. Characterizing those hyperplanes in the description of $P_M$ given by Theorem \ref{H-description-LPM-polytopes}, that either correspond to a connected contraction or deletion of $M$ on one element less or to a submatroid with one more connected component yields the following: 

\begin{corollary}\label{cor:facetdefining}
Let $M=M[U,L]$ be a connected LPM. Then the facet defining hyperplanes of $P_M$ are of the form:
\begin{itemize}
    \item[(a)] $\sum_{j=1}^{i}x_j = \sum_{j=1}^{i}L_j$ for some $1\leq i<n$ if the $i$-th point of $L$ is concave in the diagram of $M[U,L]$, or 
    \item[(b)] $\sum_{j=1}^{i}x_j = \sum_{j=1}^{i}U_j$ for some $1\leq i<n$ if the $i$-th point of $U$ is concave in the diagram of $M[U,L]$, or
    \item[(c)] $x_i=0$  for some $1\leq i\leq n$ unless the only vertical $i$-th segments among the paths in the diagram of $M[U,L]$ are the ones of $U$ and $L$, or
    \item[(d)] $x_i=1$  for some $1\leq i\leq n$  unless the only horizontal $i$-th segments among the paths in the diagram of $M[U,L]$ are the ones of $U$ and $L$.
\end{itemize}
\end{corollary}

\begin{example}
    Let us illustrate Corollary~\ref{cor:facetdefining}. Consider the matroid $M=M[12467,45678]$ depicted in~\cref{fig:concave}. Using the notation above, the point $p=(1,2)$ being concave gives rise to the facet-defining hyperplane $H_p:x_1+x_2+x_3 = 2$. The matroid corresponding to the facet $P_M\cap H_p$ is depicted in Figure \ref{fig:concave} as well as the matroid $M_q$ corresponding to the face $P_M\cap H_{(1,3)}$, where $H_{(1,3)}:x_1+x_2+x_3+x_4 = 3$. Note that $H_{(1,3)}$, despite being a defining hyperplane, is not a facet-defining one.

\begin{figure}[h]
    \centering
    \begin{tikzpicture}
             \draw[step=0.5cm, color=gray] (0, 0) grid (1,1);
             \draw[step=0.5cm, color=gray] (1,0) grid (1.5,2.5);
             \draw[step=0.5cm, color=gray] (0.5, 1) grid (1,1.5);
             \draw[line width=1pt, color=gray] (0,0) -- (0,1) -- (0.5,1) -- (0.5,1.5) -- (1,1.5) -- (1,2.5) -- (1.5,2.5);
            \draw[line width=1pt, color=gray] (0,0) -- (1.5,0) -- (1.5,2.5);
            \filldraw [color = red] (0.5,1) circle (2pt);
            \filldraw [color = blue] (0.5,1.5) circle (2pt);
    \end{tikzpicture}
    \hspace{2cm}
    \begin{tikzpicture}
             \draw[step=0.5cm, color=gray] (0, 0) grid (0.5,1);
             \draw[step=0.5cm, color=gray] (0.5, 1) grid (1,1.5);
             \draw[step=0.5cm, color=gray] (1,1) grid (1.5,2.5);
            \draw[line width=1pt,color=red] (0,0) -- (0,1) -- (0.5,1) -- (0.5,1.5) -- (1,1.5) -- (1,2.5) -- (1.5,2.5);
            \draw[line width=1pt,color=red] (0,0) -- (0.5,0) -- (0.5,1) -- (1.5,1) -- (1.5,2.5);
    \end{tikzpicture}
    \hspace{2cm}
    \begin{tikzpicture}
             \draw[step=0.5cm, color=gray] (0, 0) grid (0.5,1);
             \draw[step=0.5cm, color=gray] (1,1.5) grid (1.5,2.5);
            \draw[line width=1pt,color=blue] (0,0) -- (0,1) -- (0.5,1) -- (0.5,1.5) -- (1,1.5) -- (1,2.5) -- (1.5,2.5);
            \draw[line width=1pt,color=blue] (0,0) -- (0.5,0) -- (0.5,1.5) -- (1.5,1.5) -- (1.5,2.5);
    \end{tikzpicture}
        \caption{We show the following LPMs: $M=M[12467,45678]$ on the left, $M_{(1,2)}$ in the center, and $M_{(1,3)}$ on the right.}
    \label{fig:concave}
\end{figure}
\end{example}

\subsection{Integer and interior points of LPMs}

Let us recall some results from the literature regarding points in $P_M$ where $M$ is an LPM over $[n]$, thus $P_M\subseteq\mathbb R^n$. Given such $P_M$ we say that $Q$ is a \emph{generalized lattice path in $M$} if $Q=(s_1,\dots,s_n)$ is a sequence of $n$ line segments from $(0,0)$ to $(n-k,k)$ inside the diagram of $M$, such that the segment $s_i=\overline{AB}$ satisfies that $A=(x_{i-1},y_{i-1})$ is a point on the line $l_i=x+y=i-1$, $B=(x_{i},y_{i})\in l_{i+1}$ and $x_{i-1} \leq x_{i}$, $y_{i-1} \leq y_{i}$, for each $i\in[n]$. 
We refer to the tuple $q:=(q_1,\ldots,q_n)$ as the \emph{state vector} of
$Q$, where $q_i:=y_i-y_{i-1}$ is the vertical increment of the
segment $s_i$. Thus, the displacement vector of $s_i$ is
$(1-q_i,q_i)$.
We will refer to the points $(x_i,y_i)$ as \emph{bending points of $Q$}. Notice that $ q$ determines $Q$ and thus we use them interchangeably. 

In \cite{knauer2018lattice} the authors prove that $ q=(q_1,\dots,q_n)\in\mathbb R^n$ is a point in $P_M$ if and only if $ q$ is the state vector of a generalized lattice path in $M$. More concretely, they showed the following.
\begin{theorem}[{\cite[Theorem 3.4]{knauer2018lattice}}]\label{generalized-lattice-paths}
Let $M = M[U,L]$ be an LPM and $\mathcal{C}_M$ the set of state vectors of generalized lattice paths inside of the diagram of $M$. Then $P_M = \mathcal{C}_M.$
\end{theorem}

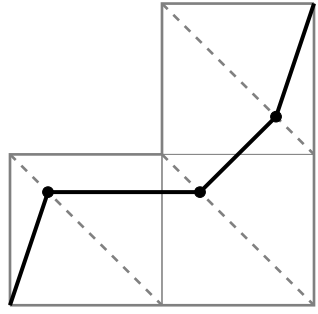
\begin{figure}[h]
    \centering
    \begin{tikzpicture}
             \draw[step=2cm, color=gray] (0, 0) grid (2,2);
             \draw[step=2cm, color=gray] (2, 0) grid (4,4);
            \draw[line width=1pt, color=gray] (0,0) -- (4,0) -- (4,4);
            \draw[line width=1pt, color=gray] (0,0) -- (0,2) -- (2,2) -- (2,4) -- (4,4);
            \draw[line width = 1pt, dashed, color=gray] (0,2) -- (2,0);
            \draw[line width = 1pt, dashed, color=gray] (2,2) -- (4,0);
            \draw[line width = 1pt, dashed, color=gray] (2,4) -- (4,2);
            \draw[line width = 1.5pt, color=black] (0,0) -- (0.5,1.5) -- (2.5,1.5) -- (3.5,2.5) -- (4,4);
            \filldraw [color = black] (0.5,1.5) circle (2pt);
            \filldraw [color = black] (2.5,1.5) circle (2pt);
            \filldraw [color = black] (3.5,2.5) circle (2pt);
            
    \end{tikzpicture}
    \caption{A generalized lattice path in $M[13,34]$.}
    \label{fig:gen-lattice-paths}
\end{figure}

For an illustration consider Figure \ref{fig:gen-lattice-paths},  showing a generalized lattice path inside the diagram of $M=M[13,34]$, whose state vector is $\left(\frac{3}{4},0,\frac{1}{2},\frac{3}{4}\right)$.

\begin{remark}\label{rem:interior_points}
Since every point of $P_M$ is the state vector of a generalized
lattice path in the diagram of $M$, a point $p\in\mathbb R^n$ belongs
to $tP_M\cap\mathbb Z^n$ if and only if $p=tq$ for the state vector
$q$ of some generalized lattice path $Q$ in $M$. Equivalently, the
bending points $(x_i,y_i)$ of $Q$ satisfy
$(tx_i,ty_i)\in\mathbb Z^2$ for every $i$.

Moreover, the description of $P_M$ in Theorem~\ref{generalized-lattice-paths} implies that a
point $q\in P_M$ belongs to $\operatorname{relint}(P_M)$ if and only
if the corresponding generalized lattice path $Q$ intersects $U$ and
$L$ only at $(0,0)$ and $(n-k,k)$ and its bending points satisfy
$x_{i-1}<x_i$ and $y_{i-1}<y_i$ for every $i\in[n]$. In other words,
$q\in\operatorname{relint}(P_M)$ if and only if $Q$ is strictly
monotone and lies strictly between $U$ and $L$, except at its
endpoints.
 \end{remark}

\subsection{Characterizing Gorenstein LPMs}
\label{subsec:characterizing-gorenstein-lpms}
We use the standard facet characterization of reflexive polytopes.
If a facet equation is written using a primitive integral normal,
then the absolute value of its constant term is the lattice distance
of the facet from the origin. A lattice polytope containing the origin
in its relative interior is reflexive if and only if all its facets
have lattice distance one; see \cite[Sections~1.1 and~1.3]{KasprzykNill}.
For a non-full-dimensional polytope, all lattice notions are taken in
the affine lattice spanned by the polytope; see
\cite[Definition~2.6 and Proposition~2.7]{LM22}.

The following proposition records exactly what this criterion means
for the facets described in Corollary~\ref{cor:facetdefining}.

\begin{proposition}
\label{prop:lpm-gorenstein-criterion}
Let $M=M[U,L]$ be a connected LPM, $\delta$ a positive
integer, and 
$z=(z_1,\ldots,z_n)\in\operatorname{relint}(\delta P_M)\cap
\mathbb Z^n$. Put $s_i:=\sum_{j=1}^i z_j$,
$\lambda_i:=\sum_{j=1}^iL_j$, and
$\mu_i:=\sum_{j=1}^iU_j$. Then $\delta P_M-z$ is reflexive if and
only if the following conditions hold:
\begin{enumerate}
    \item[(a)] $s_i-\delta\lambda_i=1$ whenever the $i$-th point of
    $L$ is concave;
    \item[(b)] $\delta\mu_i-s_i=1$ whenever the $i$-th point of $U$
    is concave;
    \item[(c)] $z_i=1$ whenever $x_i=0$ is a facet-defining
    hyperplane of $P_M$;
    \item[(d)] $\delta-z_i=1$ whenever $x_i=1$ is a facet-defining
    hyperplane of $P_M$.
\end{enumerate}
\end{proposition}

\begin{proof}
The lattice parallel to $\operatorname{aff}(P_M)$ is
$\Lambda:=\{x\in\mathbb Z^n:\sum_{j=1}^n x_j=0\}$. 
The restrictions to $\Lambda$ of the facet functionals appearing in Corollary~\ref{cor:facetdefining} are primitive elements of the dual lattice $\operatorname{Hom}(\Lambda,\mathbb Z)$.  Indeed, $x_i(e_i-e_b)=1$ for any $b\neq i$, while
$(\sum_{j=1}^ix_j)(e_a-e_b)=1$ whenever $a\leq i<b$.

By Corollary~\ref{cor:facetdefining}, the facets of $P_M$ are precisely
the relevant hyperplanes of the forms
$\sum_{j=1}^i x_j=\lambda_i$,
$\sum_{j=1}^i x_j=\mu_i$, $x_i=0$, and $x_i=1$. After dilating by
$\delta$ and translating by $z$, their lattice distances from the
origin are, respectively,
$s_i-\delta\lambda_i$, $\delta\mu_i-s_i$, $z_i$, and
$\delta-z_i$. The result follows from the facet characterization of
reflexivity recalled above.
\end{proof}

See Figure~\ref{fig:Gorenstein_LPMs} for an illustration of the following conditions on the diagrams of LPMs.
\Gorenstein*
\begin{proof}
Suppose first that $M$ is $\delta$-Gorenstein, and choose
$z\in\operatorname{relint}(\delta P_M)\cap\mathbb Z^n$ such that
$\delta P_M-z$ is reflexive.
Then $z$ satisfies 
Proposition~\ref{prop:lpm-gorenstein-criterion}. Put $q:=z/\delta$,
and let $\widetilde q$ be the generalized lattice path corresponding
to $q$ by Theorem~\ref{generalized-lattice-paths}. By
Remark~\ref{rem:interior_points}, the path $\widetilde q$ is
strictly monotone and lies strictly inside the diagram except at its
endpoints.

Suppose that the diagram of $M$ contains a $2\times2$ square. The
internal horizontal and vertical segments of this square have the
same label, say $i$. By Corollary~\ref{cor:facetdefining}, both $x_i=0$
and $x_i=1$ are facets of $P_M$. Conditions (c) and (d) of
Proposition~\ref{prop:lpm-gorenstein-criterion} give $z_i=1$ and
$\delta-z_i=1$, and hence $\delta=2$. Consequently, if
$\delta\geq3$, then the diagram contains no $2\times2$ square, so $M$
is a snake by Definition~\ref{def:snakes}.

Assume now that $\delta=2$. Since $M$ is connected, it has no loops
or coloops. Hence every point in $\operatorname{relint}(P_M)$ has all
coordinates strictly between $0$ and $1$. It follows that $0<z_i<2$
for every $i$, and therefore $z=(1,\ldots,1)$. Since every point of
$2P_M$ has coordinate sum $2k$, we obtain $n=2k$.

The path corresponding to $z/2$ has north increment $1/2$ at every
step, and its $i$-th bending point is $(i/2,i/2)$. Thus it is the
diagonal segment $\ell$. Since $z/2\in\operatorname{relint}(P_M)$,
Remark~\ref{rem:interior_points} shows that $\ell$ intersects
$U$ and $L$ only at its endpoints.

Let $(i-\lambda_i,\lambda_i)$ be a concave point of $L$. Since
$s_i=i$, condition (a) of
Proposition~\ref{prop:lpm-gorenstein-criterion} gives
$i-2\lambda_i=1$, so this point is of the form $(r+1,r)$. Similarly,
condition (b) shows that every concave point of $U$ is of the form
$(r,r+1)$. This proves necessity when $\delta=2$.

Conversely, suppose that the conditions in the case $\delta=2$ hold,
and set $z=(1,\ldots,1)$. The point $z/2$ corresponds to the diagonal
path, and the intersection assumption together with
Remark~\ref{rem:interior_points} gives
$z\in\operatorname{relint}(2P_M)$. Conditions (c) and (d) of
Proposition~\ref{prop:lpm-gorenstein-criterion} hold automatically.
If $(i-\lambda_i,\lambda_i)=(r+1,r)$ is a concave point of $L$, then
$i=2r+1$ and $i-2\lambda_i=1$, so condition (a) holds. The same
calculation for a concave point $(r,r+1)$ of $U$ gives condition (b).
Thus Proposition~\ref{prop:lpm-gorenstein-criterion} shows that
$2P_M-z$ is reflexive, and $M$ is $2$-Gorenstein.

It remains to consider $\delta\geq3$. We have already shown that $M$
is a snake. Up to reflection across $\ell$, write
$M=S(\alpha_1,\ldots,\alpha_s)$ using the notation established below 
Definition~\ref{def:snakes}. 
The labels split into
consecutive blocks $I_1,\ldots,I_s$ corresponding to the successive
strips. If $s\geq2$, then $|I_1|=\alpha_1+1$,
$|I_j|=\alpha_j$ for $2\leq j\leq s-1$, and
$|I_s|=\alpha_s+1$. If $s=1$, then $|I_1|=\alpha_1+2$.

For a label $i$ in a horizontal strip, $x_i=0$ is a facet, so
condition (c) gives $q_i=1/\delta$. For a label $i$ in a vertical
strip, $x_i=1$ is a facet, so condition (d) gives
$q_i=(\delta-1)/\delta$. Thus each step of $\widetilde q$ advances by
$1/\delta$ in the short direction of its strip. Conditions (a) and
(b) say that $\widetilde q$ passes at distance $1/\delta$ from every
concave corner.

Suppose first that $s\geq2$. In the first strip, the displacement of
$\widetilde q$ in the short direction is both
$(\alpha_1+1)/\delta$ and $1-1/\delta$, so
$\alpha_1=\delta-2$. In an internal strip, it is both
$\alpha_j/\delta$ and $1-2/\delta$, so
$\alpha_j=\delta-2$. The same calculation in the last strip gives
$\alpha_s=\delta-2$. If $s=1$, the path crosses the short direction
of the unique strip completely, so
$(\alpha_1+2)/\delta=1$ and again $\alpha_1=\delta-2$. Hence
$M\cong S(\delta-2,\ldots,\delta-2)$, or, after reflection,
$M\cong S^*(\delta-2,\ldots,\delta-2)$.

Conversely, let $M=S(\delta-2,\ldots,\delta-2)$ for some
$\delta\geq3$. Define a generalized path $\widetilde q$ by taking
$q_i=1/\delta$ on horizontal strips and
$q_i=(\delta-1)/\delta$ on vertical strips. The preceding strip
calculations, read in reverse, show that $\widetilde q$ is strictly
monotone, lies strictly inside the diagram between its endpoints, and
passes at distance $1/\delta$ from every concave corner. Hence
$z:=\delta q$ belongs to
$\operatorname{relint}(\delta P_M)\cap\mathbb Z^n$ by
Theorem~\ref{generalized-lattice-paths} and
Remark~\ref{rem:interior_points}. Conditions (a)--(d) of
Proposition~\ref{prop:lpm-gorenstein-criterion} now hold by
construction. Therefore $\delta P_M-z$ is reflexive. The case
$S^*(\delta-2,\ldots,\delta-2)$ follows by reflection.
\end{proof}

\begin{figure}[h]
    \centering

    \begin{tikzpicture}
        \node at (0,0) {\begin{tikzpicture}
            \draw[line width=1pt, color=gray] (0,0) -- (0,0.75) -- (1.5,0.75) -- (1.5,2.25) -- (3.75,2.25);
            \draw[line width=1pt, color=gray] (0,0) -- (2.25,0) -- (2.25,1.5) -- (3.75,1.5)-- (3.75,2.25);
            
            \draw[step=0.75cm, color=gray] (0, 0) grid (2.25,0.75);
            \draw[step=0.75cm, color=gray] (1.5, 0) grid (2.25,2.25);
            \draw[step=0.75cm, color=gray] (1.5, 1.5) grid (3.75,2.25);
       \end{tikzpicture}};
    
        \node at (3.75,0) {\begin{tikzpicture}
                \draw[line width=1pt, color=gray] (0,0) -- (0,1.5) -- (0.75,1.5) -- (0.75,3) -- (3,3);
                \draw[line width=1pt, color=gray] (0,0) -- (1.5,0) -- (1.5,0.75) -- (2.25,0.75) -- (2.25,1.5) -- (3,1.5) -- (3,3);
                
                \draw[step=0.75cm, color=gray] (0, 0) grid (0.75,1.5);
                \draw[step=0.75cm, color=gray] (0.75, 0) grid (1.5,3);
                \draw[step=0.75cm, color=gray] (1.5, 0.75) grid (2.25,3);
                \draw[step=0.75cm, color=gray] (2.25,1.5) grid (3,3);
               \draw[line width=1pt, dashed, color = red] (-0.25,-0.25) -- (3.25,3.25);
               \node[red] at (-0.45,-0.45) {$\ell$};
       \end{tikzpicture}};
    
        \node at (9.25,0) {\begin{tikzpicture}
                \draw[step=0.75cm, color=gray] (0, 0) grid (1.5,2.25);
                \draw[step=0.75cm, color=gray] (1.5, 0.75) grid (2.25,3);
                \draw[step=0.75cm, color=gray] (2.25,1.5) grid (3,3);
               \draw[line width=1pt, color = gray] (0,0) -- (0,2.25) -- (0.75,2.25) -- (0.75,3) -- (3,3);
               \draw[line width=1pt, color = gray] (0,0) -- (1.5,0) -- (1.5,0.75) -- (2.25,0.75) -- (2.25,1.5) -- (3,1.5) -- (3,3);
               \draw[line width=1pt, dashed, color = red] (-0.25,-0.25) -- (3.25,3.25);
               \node[red] at (-0.45,-0.45) {$\ell$};
               \filldraw [color = black] (0.75,2.25) circle (2pt);
       \end{tikzpicture}};
    
        \node at (12.75,0) {\begin{tikzpicture}
                \draw[step=0.75cm, color=gray] (0.75, 0) grid (1.5,3);
                \draw[step=0.75cm, color=gray] (1.5, 0.75) grid (2.25,3);
                \draw[step=0.75cm, color=gray] (2.25,1.5) grid (3,3);
               \draw[line width=1pt, color = gray] (0,0) -- (0,0.75) -- (0.75,0.75) -- (0.75,3) -- (3,3);
               \draw[line width=1pt, color = gray] (0,0) -- (1.5,0) -- (1.5,0.75) -- (2.25,0.75) -- (2.25,1.5) -- (3,1.5) -- (3,3);
               
               \draw[line width=1pt, dashed, color = red] (-0.25,-0.25) -- (3.25,3.25);
               \node[red] at (-0.45,-0.45) {$\ell$};
               \filldraw [color = black] (0.75,0.75) circle (2pt);
       \end{tikzpicture}};
    \end{tikzpicture}
   
    \caption{The two LPMs on the left are Gorenstein according to Theorem~\ref{thm:Gorenstein-LPM}. The leftmost one is $S^*(3,3,3)$, so $\delta=5$, and the concave points of the other one are $\{(1,2),(2,1),(3,2)\}$. The two LPMs on the right are not Gorenstein as the highlighted concave points show.}
    \label{fig:Gorenstein_LPMs}
\end{figure}

As uniform matroids are LPMs, we recover a result of \cite[Theorem 2.4]{DENEGRI1997629}.

\begin{corollary}
A uniform matroid is Gorenstein if and only if it is $U_{n,2n}$, $U_{1,n}$ or $U_{n-1,n}$. 
\end{corollary}

Finally, and as a main connection to the theme of the paper we obtain:

\begin{corollary}
If an LPM satisfies the conditions of Theorem~\ref{thm:Gorenstein-LPM}, then its $h^*$-vector is unimodal.
\end{corollary}
\begin{proof}
Let $M$ be an LPM. Since $P_M$ is alcoved it is folklore that the alcoved triangulation of $P_M$ is unimodular and regular. Hence, one can apply~\cite[Theorem 1]{bruns2007h}, which guarantees that Gorenstein polytopes with such triangulations have unimodal $h^*$-vector.
\end{proof}

To finish this section, we recall that positroid polytopes are the intersection of matroid polytopes and alcoved polytopes~\cite[Theorem 2.1]{lam2020polypositroids}. We are interested in extending our results to the family of positroids.
More specifically we are interested in the following problem that would be a direct generalization of \cref{thm:Gorenstein-LPM}.

\begin{question}\label{ques:PositroidGorenstein}
    Characterize all positroids that are Gorenstein. 
\end{question}


\section{Towards positroidal snakes}\label{sec:positroids}

The motivation for the following discussion complements the role of snakes provided by \cref{thm:snakesposets,thm:Gorenstein-LPM}.
As mentioned earlier, the importance of snakes mainly stems from the fact that the matroid polytope $P_M$ of an LPM $M$ decomposes into snakes~\cite{chatelain2011matroid}. More precisely, let $\Sigma_{\mathrm{box}}(M)$ denote the
subdivision of $P_M$ induced by the hyperplanes $y_i=\sum_{j=1}^ix_j\in\mathbb Z$. One direction of~\cref{prop:order-matroid-subdivision} formalizes the  fact that
\begin{itemize}
    \item[(i)] the map $x\mapsto (\sum_{j=1}^ix_j)_{1\leq i<n}$ sends pieces of $\Sigma_{\mathrm{box}}(M)$ to integer translates of order polytopes,
    \item[(ii)] $\Sigma_{\mathrm{box}}(M)$ is a (finest) 
    decomposition of $P_M$ into matroid polytopes.
\end{itemize}

Since every integer translate of an order polytope is a union of alcoves, (i) is equivalent to $P_M$ being an alcoved polytope by \cite[Section~2.3]{Lam2007}. Hence, (i) is equivalent to $M$ being a positroid.

We now show that (ii) is equivalent to $M$ being an LPM.

\begin{lemma}\label{lem:snake-matroid-indecomposable}
Let $S$ and $N$ be matroids on $[n]$.  If $S$ is a snake, $P_N\subseteq P_S$ and $
    \dim P_N=\dim P_S=n-1$, then $N=S$. 
\end{lemma}

\begin{proof}
Clearly,
$P_N\subseteq P_S$ gives $\mathcal B(N)\subseteq\mathcal B(S)$.
Thus, the identity on $[n]$ is a rank-preserving weak map from $S$ to $N$.
Moreover, a matroid on $[n]$ with $q$ connected components has
base-polytope dimension $n-q$, so $S$ and $N$ are connected. Since $S$ is a snake, it is graphic and hence binary
\cite[Theorem~3.2 and Corollary~3.3]{KMR18}. 
By~\cite[Theorem~6.10]{Lucas1975}, a nontrivial
rank-preserving weak image of a binary matroid is disconnected. Hence the
weak map is trivial and $N=S$.
\end{proof}

\OrderMatroidSubdivision*

\begin{proof}
The first direction was already known. Indeed, suppose first that $M=M[U,L]$ is a lattice path matroid of rank $k$. By
\cite[Theorem~1 and Corollary~4]{chatelain2011matroid} $P_M$ can be split into snake polytopes, see
also \cite[Theorem~2.6]{FMP26}. Indeed, the hyperplanes occurring in the LPM splits are exactly the hyperplanes $y_i=\sum_{j=1}^ix_j\in\mathbb Z$ that cut the relative interior of $P_M$. Hence the resulting snake
subdivision is $\Sigma_{\mathrm{box}}(M)$ and matroidal in particular.

Conversely, suppose that $\Sigma_{\mathrm{box}}(M)$ is matroidal.
Since the uniform matroid $U_{k,n}$ is a lattice path matroid, the first
part of the proof applied to $U_{k,n}$ shows that the maximal cells of
$\Sigma_{\mathrm{box}}(U_{k,n})$ are snake matroid polytopes.

Let $P_{S_a}$ be such a maximal cell for which
 $P_N=P_M\cap P_{S_a}$ 
is a maximal cell of $\Sigma_{\mathrm{box}}(M)$. Since $M$ is connected,
$\dim P_M=n-1$, and therefore
$
        \dim P_N=\dim P_{S_a}=n-1$. 
As $P_N\subseteq P_{S_a}$,
\cref{lem:snake-matroid-indecomposable} gives $N=S_a$.
Consequently, every full-dimensional cell of
$\Sigma_{\mathrm{box}}(U_{k,n})$ met by
$\operatorname{relint}(P_M)$ is contained in $P_M$ in its entirety.
Thus $P_M$ is a union of whole cells of
$\Sigma_{\mathrm{box}}(U_{k,n})$.

Let $F$ be a facet of $P_M$. A generic point of
$\operatorname{relint}(F)$ lies in the relative interior of a facet of
one of these ambient cells. Hence the supporting hyperplane of $F$ is
either inherited from $\Delta(k,n)$ or has the form
$\sum_{j=1}^ix_j=m$ 
for some $1\leq i\leq n$ and $m\in\mathbb Z$. 
Defining $\ell_i=\min_{x\in P_M} \sum_{j=1}^ix_j$ and $        u_i=\max_{x\in P_M} \sum_{j=1}^ix_j$ for $1\leq i<n$, it follows that
\[
 P_M=
 \left\{x\in\Delta(k,n):
        \ell_i\leq \sum_{j=1}^ix_j\leq u_i
        \text{ for }1\leq i<n\right\}.                 \tag{$\star$}
\]
Set $\ell_0=u_0=0$ and $\ell_n=u_n=k$. Since the vertices of $P_M$
are the indicator vectors of the bases of $M$,
 $u_i=r_M([i])$ and $\ell_i=k-r_M(\{i+1,\ldots,n\})$. 
Hence the successive differences of both $(\ell_i)_{i=0}^n$ and
$(u_i)_{i=0}^n$ belong to $\{0,1\}$. They therefore define lattice
paths $L$ and $U$, with $U$ weakly above $L$. By
\cref{H-description-LPM-polytopes}, the right-hand side of~$(\star)$ is
$P_{M[U,L]}$. Thus $P_M=P_{M[U,L]}$
and comparison of the $0/1$ vertices yields $M=M[U,L]$.
Therefore $M$ is a lattice path matroid.

Finally, in the lattice path case the maximal cells of
$\Sigma_{\mathrm{box}}(M)$ are snakes, and
\cref{lem:snake-matroid-indecomposable} shows that none of them admits
a nontrivial matroid subdivision. Hence
$\Sigma_{\mathrm{box}}(M)$ is a finest matroid subdivision.
\end{proof}

\begin{figure}[ht]
    \centering
    \includegraphics[width=.55\textwidth]{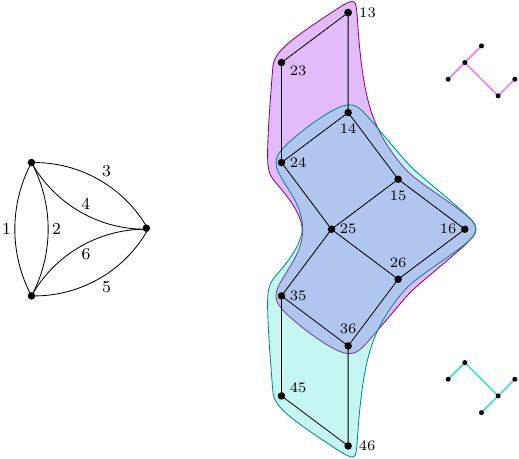}
    \caption{A graph representing a positroid $M$ that is not an LPM, its set of bases decomposed via $\Sigma_{\mathrm{box}}(M)$, and the associated posets.}\label{fig:nonLPMpoistroid}
\end{figure}

We illustrate $\Sigma_{\mathrm{box}}(M)$ for the following positroid in Figure~\ref{fig:nonLPMpoistroid}. Let $M$ be the graphic matroid of the triangle with three pairs of double edges $e_1,e_2$ and $e_3,e_4$ and $e_5,e_6$. Then $M$ is a positroid, see~\cite[Theorem 5.1]{Blum2001BasesortableMA}. 
However, $M$ is not an LPM, see~\cite[Theorem 3.2]{KMR18}. By Theorem~\ref{thm:snakesposets} $P_M$ is not an order polytope and $\Sigma_{\mathrm{box}}(M)$ consists of two order polytopes, neither of which is a matroid polytope since their underlying poset is not a fence. On the other hand $P_M$ is indecomposable into smaller matroid polytopes, see~\cite[Example 7.9]{BJR09}. Based on this example, we formulate a final question related to coarsest subdivisions into order polytopes:

\begin{question}
    Which posets arise in the subdivision $\Sigma_{\mathrm{box}}(M)$ of positroids into order polytopes?
\end{question}

\subsubsection*{Acknowledgements} We  thank Vic Reiner and Volkmar Welker for explaining to us their results in~\cite{RW05}. C.~Benedetti thanks grant FAPA of Universidad de los Andes.
K.~Knauer was supported through grant 
PID2022-137283NB-C22 funded by MICIU/AEI/10.13039/501100011033 by ERDF/EU and through the Severo Ochoa and Mar\'ia de Maeztu Program for Centers and Units of Excellence in R\&D (CEX2020-001084-M). J. Valencia-Porras was partially funded by the Natural Sciences and Engineering Research Council of Canada (NSERC) through Discovery Grants RGPIN-2026-04714 and RGPIN-2021-02568.

\printbibliography

@article{An17,
  author        = {An, Suhyung and Jung, JiYoon and Kim, Sangwook},
  title         = {Facial structures of lattice path matroid polytopes},
  journal       = {Discrete Math.},
  volume        = {343},
  number        = {1},
  pages         = {111628},
  year          = {2020},
  doi           = {10.1016/j.disc.2019.111628}
}

@article{ardila2010matroid,
  author        = {Ardila, Federico and Benedetti, Carolina and Doker, Jeffrey},
  title         = {Matroid polytopes and their volumes},
  journal       = {Discrete Comput. Geom.},
  volume        = {43},
  number        = {4},
  pages         = {841--854},
  year          = {2010},
  doi           = {10.1007/s00454-009-9232-9}
}

@incollection{BatyrevNill,
  author        = {Batyrev, Victor and Nill, Benjamin},
  title         = {Combinatorial aspects of mirror symmetry},
  booktitle     = {Integer Points in Polyhedra---Geometry, Number Theory, Representation Theory, Algebra, Optimization, Statistics},
  series        = {Contemporary Mathematics},
  volume        = {452},
  publisher     = {American Mathematical Society},
  address       = {Providence, RI},
  year          = {2008},
  pages         = {35--66},
  doi           = {10.1090/conm/452/08770}
}

@article{BK22,
  author        = {Benedetti-Vel{\'{a}}squez, Carolina and Knauer, Kolja},
  title         = {Lattice path matroids and quotients},
  journal       = {Combinatorica},
  volume        = {44},
  number        = {3},
  pages         = {621--650},
  year          = {2024},
  doi           = {10.1007/s00493-024-00085-4}
}

@misc{Bid-12,
  author        = {Bidkhori, Hoda},
  title         = {Lattice Path Matroid Polytopes},
  year          = {2012},
  eprint        = {1212.5705},
  archivePrefix = {arXiv},
  primaryClass  = {math.CO}
}

@article{BJR09,
  author        = {Billera, Louis J. and Jia, Ning and Reiner, Victor},
  title         = {A quasisymmetric function for matroids},
  journal       = {Eur. J. Combin.},
  volume        = {30},
  number        = {8},
  pages         = {1727--1757},
  year          = {2009},
  doi           = {10.1016/j.ejc.2008.12.007}
}

@article{Blum2001BasesortableMA,
  author        = {Blum, Stefan},
  title         = {Base-sortable matroids and {K}oszulness of semigroup rings},
  journal       = {Eur. J. Combin.},
  volume        = {22},
  number        = {7},
  pages         = {937--951},
  year          = {2001},
  doi           = {10.1006/eujc.2001.0526}
}

@article{BONIN2006701,
  author        = {Bonin, Joseph E. and de Mier, Anna},
  title         = {Lattice path matroids: structural properties},
  journal       = {Eur. J. Combin.},
  volume        = {27},
  number        = {5},
  pages         = {701--738},
  year          = {2006},
  doi           = {10.1016/j.ejc.2005.01.008}
}

@article{bonin2003lattice,
  author        = {Bonin, Joseph E. and de Mier, Anna and Noy, Marc},
  title         = {Lattice path matroids: {E}numerative aspects and {T}utte polynomials},
  journal       = {J. Combin. Theory Ser. A},
  volume        = {104},
  number        = {1},
  pages         = {63--94},
  year          = {2003},
  doi           = {10.1016/S0097-3165(03)00122-5}
}

@article{Bon-10,
  author        = {Bonin, Joseph E.},
  title         = {Lattice path matroids: the excluded minors},
  journal       = {J. Combin. Theory Ser. B},
  volume        = {100},
  number        = {6},
  pages         = {585--599},
  year          = {2010},
  doi           = {10.1016/j.jctb.2010.05.001}
}

@article{B04,
  author        = {Br{\"{a}}nd{\'{e}}n, Petter},
  title         = {Counterexamples to the {N}eggers--{S}tanley conjecture},
  journal       = {Electron. Res. Announc. Amer. Math. Soc.},
  volume        = {10},
  pages         = {155--158},
  year          = {2004},
  doi           = {10.1090/S1079-6762-04-00140-4}
}

@article{bruns2007h,
  author        = {Bruns, Winfried and R{\"{o}}mer, Tim},
  title         = {$h$-vectors of {G}orenstein polytopes},
  journal       = {J. Combin. Theory Ser. A},
  volume        = {114},
  number        = {1},
  pages         = {65--76},
  year          = {2007},
  doi           = {10.1016/j.jcta.2006.03.003}
}

@article{chatelain2011matroid,
  author        = {Chatelain, Vanessa and Ram{\'{\i}}rez Alfons{\'{\i}}n, Jorge Luis},
  title         = {Matroid base polytope decomposition},
  journal       = {Adv. Appl. Math.},
  volume        = {47},
  number        = {1},
  pages         = {158--172},
  year          = {2011},
  doi           = {10.1016/j.aam.2010.04.005}
}

@article{DeLoera2008,
  author        = {De Loera, Jes{\'{u}}s A. and Haws, David C. and K{\"{o}}ppe, Matthias},
  title         = {Ehrhart polynomials of matroid polytopes and polymatroids},
  journal       = {Discrete Comput. Geom.},
  volume        = {42},
  number        = {4},
  pages         = {670--702},
  year          = {2009},
  doi           = {10.1007/s00454-008-9080-z}
}

@article{DeM-07,
  author        = {De Mier, Anna},
  title         = {A natural family of flag matroids},
  journal       = {SIAM J. Discrete Math.},
  volume        = {21},
  number        = {1},
  pages         = {130--140},
  year          = {2007},
  doi           = {10.1137/050638515}
}

@article{DENEGRI1997629,
  author        = {De Negri, Emanuela and Hibi, Takayuki},
  title         = {Gorenstein algebras of {V}eronese type},
  journal       = {J. Algebra},
  volume        = {193},
  number        = {2},
  pages         = {629--639},
  year          = {1997},
  doi           = {10.1006/jabr.1997.6990}
}

@article{Del-12,
  author        = {Delucchi, Emanuele and Dlugosch, Martin},
  title         = {Bergman complexes of lattice path matroids},
  journal       = {SIAM J. Discrete Math.},
  volume        = {29},
  number        = {4},
  pages         = {1916--1930},
  year          = {2015},
  doi           = {10.1137/130944242}
}

@article{Ehrhart,
  author        = {Ehrhart, Eug{\`e}ne},
  title         = {Sur les poly{\`e}dres rationnels homoth{\`e}tiques {\`a} $n$ dimensions},
  journal       = {C. R. Acad. Sci. Paris},
  volume        = {254},
  pages         = {616--618},
  year          = {1962}
}

@misc{FMP26,
  author        = {Ferroni, Luis and Morales, Alejandro H. and Panova, Greta},
  title         = {Ehrhart positivity for lattice path matroids},
  year          = {2026},
  eprint        = {2605.22673},
  archivePrefix = {arXiv},
  primaryClass  = {math.CO}
}

@article{ferroni2021ehrhart,
  author        = {Ferroni, Luis and Jochemko, Katharina and Schr{\"{o}}ter, Benjamin},
  title         = {Ehrhart polynomials of rank two matroids},
  journal       = {Adv. Appl. Math.},
  volume        = {141},
  pages         = {102410},
  year          = {2022},
  doi           = {10.1016/j.aam.2022.102410}
}

@article{Ferroni_minimal,
  author        = {Ferroni, Luis},
  title         = {On the {E}hrhart polynomial of minimal matroids},
  journal       = {Discrete Comput. Geom.},
  volume        = {68},
  number        = {1},
  pages         = {255--273},
  year          = {2021},
  doi           = {10.1007/s00454-021-00313-4}
}

@article{GGMS,
  author        = {Gelfand, Israel M. and Goresky, R. Mark and MacPherson, Robert and Serganova, Vera},
  title         = {Combinatorial geometries, convex polyhedra and {S}chubert cells},
  journal       = {Adv. Math.},
  volume        = {63},
  number        = {3},
  pages         = {301--316},
  year          = {1987},
  doi           = {10.1016/0001-8708(87)90059-4}
}

@article{HIBI2017991,
  author        = {Hibi, Takayuki and Li, Nan and Sahara, Yoshimi and Shikama, Akihiro},
  title         = {The numbers of edges of the order polytope and the chain polytope of a finite partially ordered set},
  journal       = {Discrete Math.},
  volume        = {340},
  number        = {5},
  pages         = {991--994},
  year          = {2017},
  doi           = {10.1016/j.disc.2017.01.005}
}

@article{Hibi2021,
  author        = {Hibi, Takayuki and Laso{\'{n}}, Micha{\l} and Matsuda, Kazunori and Micha{\l}ek, Mateusz and Vodi{\v{c}}ka, Martin},
  title         = {Gorenstein graphic matroids},
  journal       = {Israel J. Math.},
  volume        = {243},
  number        = {1},
  pages         = {1--26},
  year          = {2021},
  doi           = {10.1007/s11856-021-2136-y}
}

@inproceedings{hibi1987distributive,
  author        = {Hibi, Takayuki},
  title         = {Distributive Lattices, Affine Semigroup Rings and Algebras with Straightening Laws},
  booktitle     = {Commutative Algebra and Combinatorics},
  editor        = {Nagata, Masayoshi and Matsumura, Hideyuki},
  series        = {Advanced Studies in Pure Mathematics},
  volume        = {11},
  publisher     = {Kinokuniya Company Ltd. / North-Holland},
  address       = {Tokyo / Amsterdam},
  year          = {1987},
  pages         = {93--109},
  doi           = {10.2969/aspm/01110093}
}

@article{HolzmannNortonTobey,
  author        = {Holzmann, Carlos A. and Norton, Peter G. and Tobey, Malcolm D.},
  title         = {A graphical representation of matroids},
  journal       = {SIAM J. Appl. Math.},
  volume        = {25},
  number        = {4},
  pages         = {618--627},
  year          = {1973},
  doi           = {10.1137/0125060}
}

@article{KasprzykNill,
  author        = {Kasprzyk, Alexander M. and Nill, Benjamin},
  title         = {Reflexive polytopes of higher index and the number 12},
  journal       = {Electron. J. Combin.},
  volume        = {19},
  number        = {3},
  pages         = {Paper 9, 18},
  year          = {2012},
  doi           = {10.37236/2366}
}

@article{KMR18,
  author        = {Knauer, Kolja and Mart{\'{i}}nez-Sandoval, Leonardo and Ram{\'{\i}}rez Alfons{\'{\i}}n, Jorge Luis},
  title         = {A {T}utte polynomial inequality for lattice path matroids},
  journal       = {Adv. Appl. Math.},
  volume        = {94},
  pages         = {23--38},
  year          = {2018},
  doi           = {10.1016/j.aam.2016.11.008}
}

@article{knauer2018lattice,
  author        = {Knauer, Kolja and Mart{\'{i}}nez-Sandoval, Leonardo and Ram{\'{\i}}rez Alfons{\'{\i}}n, Jorge Luis},
  title         = {On lattice path matroid polytopes: integer points and {E}hrhart polynomial},
  journal       = {Discrete Comput. Geom.},
  volume        = {60},
  number        = {3},
  pages         = {698--719},
  year          = {2018},
  doi           = {10.1007/s00454-018-9965-4}
}

@article{Klbl2020,
  author        = {K{\"{o}}lbl, Max},
  title         = {Gorenstein graphic matroids from multigraphs},
  journal       = {Ann. Comb.},
  volume        = {24},
  number        = {2},
  pages         = {395--403},
  year          = {2020},
  doi           = {10.1007/s00026-020-00495-3}
}

@article{lam2020polypositroids,
  author        = {Lam, Thomas and Postnikov, Alexander},
  title         = {Polypositroids},
  journal       = {Forum Math. Sigma},
  volume        = {12},
  pages         = {Paper No. e42, 67},
  year          = {2024},
  doi           = {10.1017/fms.2024.11}
}

@article{Lam2007,
  author        = {Lam, Thomas and Postnikov, Alexander},
  title         = {Alcoved polytopes. {I}},
  journal       = {Discrete Comput. Geom.},
  volume        = {38},
  number        = {3},
  pages         = {453--478},
  year          = {2007},
  doi           = {10.1007/s00454-006-1294-3}
}

@article{LM22,
  author        = {Laso{\'{n}}, Micha{\l} and Micha{\l}ek, Mateusz},
  title         = {Gorenstein Matroids},
  journal       = {Int. Math. Res. Not. IMRN},
  volume        = {2023},
  number        = {18},
  pages         = {15687--15728},
  year          = {2023},
  doi           = {10.1093/imrn/rnac292}
}

@article{Lucas1975,
  author        = {Lucas, Dean},
  title         = {Weak Maps of Combinatorial Geometries},
  journal       = {Trans. Amer. Math. Soc.},
  volume        = {206},
  year          = {1975},
  pages         = {247--279},
  doi           = {10.2307/1997156}
}

@misc{benedetti2023latticepathmatroidpolytopes,
  author        = {Benedetti, Carolina and Knauer, Kolja and Valencia-Porras, Jer{\'{o}}nimo},
  title         = {On lattice path matroid polytopes: alcoved triangulations and snake decompositions},
  year          = {2023},
  eprint        = {2303.10458v1},
  archivePrefix = {arXiv},
  primaryClass  = {math.CO},
  url           = {https://arxiv.org/abs/2303.10458v1}
}

@book{Maf81,
  author        = {Maffioli, Francesco},
  title         = {An introduction to matroid optimization},
  series        = {Quaderni Serie III},
  volume        = {129},
  publisher     = {Istituto per le Applicazioni del Calcolo},
  address       = {Rome},
  year          = {1981}
}

@article{MAURER1973216,
  author        = {Maurer, Stephen B.},
  title         = {Matroid basis graphs. {I}},
  journal       = {J. Combin. Theory Ser. B},
  volume        = {14},
  number        = {3},
  pages         = {216--240},
  year          = {1973},
  doi           = {10.1016/0095-8956(73)90005-1}
}

@article{MSS21,
  author        = {McConville, Thomas and Sagan, Bruce E. and Smyth, Clifford},
  title         = {On a rank-unimodality conjecture of {M}orier-{G}enoud and {O}vsienko},
  journal       = {Discrete Math.},
  volume        = {344},
  number        = {8},
  pages         = {112483},
  year          = {2021},
  doi           = {10.1016/j.disc.2021.112483}
}

@article{Mor-13,
  author        = {Morton, Jason and Turner, Jacob},
  title         = {Computing the {T}utte polynomial of lattice path matroids using determinantal circuits},
  journal       = {Theoret. Comput. Sci.},
  volume        = {598},
  pages         = {150--156},
  year          = {2015},
  doi           = {10.1016/j.tcs.2015.07.042}
}

@article{N78,
  author        = {Neggers, Joseph},
  title         = {Representations of finite partially ordered sets},
  journal       = {J. Comb. Inf. Syst. Sci.},
  volume        = {3},
  pages         = {113--133},
  year          = {1978}
}

@book{O11,
  author        = {Oxley, James G.},
  title         = {Matroid Theory},
  edition       = {2nd},
  series        = {Oxford Graduate Texts in Mathematics},
  volume        = {21},
  publisher     = {Oxford University Press},
  address       = {Oxford},
  year          = {2011},
  doi           = {10.1093/acprof:oso/9780199603398.001.0001}
}

@article{shardPPR23,
  author        = {Padrol, Arnau and Pilaud, Vincent and Ritter, Julian},
  title         = {Shard polytopes},
  journal       = {Int. Math. Res. Not. IMRN},
  volume        = {2023},
  number        = {9},
  pages         = {7686--7796},
  year          = {2023},
  doi           = {10.1093/imrn/rnac042}
}

@article{PinedaVillavicencioSchroeter,
  author        = {Pineda-Villavicencio, Guillermo and Schr{\"{o}}ter, Benjamin},
  title         = {Reconstructibility of matroid polytopes},
  journal       = {SIAM J. Discrete Math.},
  volume        = {36},
  number        = {1},
  pages         = {490--508},
  year          = {2022},
  doi           = {10.1137/21M1401176}
}

@misc{postnikov2006total,
  author        = {Postnikov, Alexander},
  title         = {Total positivity, {G}rassmannians, and networks},
  year          = {2006},
  eprint        = {math/0609764},
  archivePrefix = {arXiv},
  primaryClass  = {math.CO}
}

@article{RW05,
  author        = {Reiner, Victor and Welker, Volkmar},
  title         = {On the {C}harney--{D}avis and {N}eggers--{S}tanley conjectures},
  journal       = {J. Combin. Theory Ser. A},
  volume        = {109},
  number        = {2},
  pages         = {247--280},
  year          = {2005},
  doi           = {10.1016/j.jcta.2004.09.003}
}

@misc{sanchez2023mirkovicvilonen,
  author        = {Sanchez, Mario},
  title         = {Mirkovi{\'{c}}--{V}ilonen Polytopes from Combinatorics},
  year          = {2023},
  eprint        = {2311.16979},
  archivePrefix = {arXiv},
  primaryClass  = {math.CO}
}

@article{Sch-10,
  author        = {Schweig, Jay},
  title         = {On the $h$-vector of a lattice path matroid},
  journal       = {Electron. J. Combin.},
  volume        = {17},
  number        = {1},
  pages         = {Research Paper N3, 6},
  year          = {2010},
  doi           = {10.37236/275}
}

@article{Sch-11,
  author        = {Schweig, Jay},
  title         = {Toric ideals of lattice path matroids and polymatroids},
  journal       = {J. Pure Appl. Algebra},
  volume        = {215},
  number        = {11},
  pages         = {2660--2665},
  year          = {2011},
  doi           = {10.1016/j.jpaa.2011.03.010}
}

@incollection{Stanley-NonNegThm,
  author        = {Stanley, Richard P.},
  title         = {Decompositions of rational convex polytopes},
  booktitle     = {Combinatorial Mathematics, Optimal Designs and Their Applications},
  series        = {Annals of Discrete Mathematics},
  volume        = {6},
  publisher     = {North-Holland},
  pages         = {333--342},
  year          = {1980},
  doi           = {10.1016/S0167-5060(08)70717-9}
}

@article{Stanley1986,
  author        = {Stanley, Richard P.},
  title         = {Two poset polytopes},
  journal       = {Discrete Comput. Geom.},
  volume        = {1},
  number        = {1},
  pages         = {9--23},
  year          = {1986},
  doi           = {10.1007/BF02187680}
}

@article{S07,
  author        = {Stembridge, John R.},
  title         = {Counterexamples to the poset conjectures of {N}eggers, {S}tanley, and {S}tembridge},
  journal       = {Trans. Amer. Math. Soc.},
  volume        = {359},
  number        = {3},
  pages         = {1115--1128},
  year          = {2007},
  doi           = {10.1090/S0002-9947-06-04271-1}
}

\end{document}